\documentclass[authoryear,preprint,11pt]{elsarticle}
\usepackage{multirow}

\usepackage{color,xcolor}
\usepackage[utf8]{inputenc}
\usepackage{tikz}
\usepackage{enumitem}
\usepackage{amsmath,amsfonts,amsthm}
\usepackage{gensymb}
\usetikzlibrary{shapes.geometric, arrows, positioning, fit, shapes, arrows.meta}
\tikzstyle{bigbox} = [draw, dashed, rounded corners, text width=18em, minimum height=8em, align=right]

\tikzstyle{bigbox1} = [draw, dashed, rounded corners, text width=18em, minimum height=4em, align=right]
\tikzstyle{startend} = [rectangle, draw, fill=!40, text width=3em, text centered, rounded corners, minimum height=2em]
\tikzstyle{block} = [circle, draw, fill=gray!25, text width=7.5em, text centered, minimum height=2em, font]
\tikzstyle{block2} = [circle, draw, fill=white!25, text width=7.5em, text centered, minimum height=2em, font]
\tikzstyle{block3} = [rectangle, text width=20em, text centered, rounded corners, minimum height=2em]
\tikzstyle{block4} = [rectangle, draw, fill=gray!15, text width=6em, text centered, rounded corners, minimum height=2em]
\tikzstyle{block5} = [rectangle, draw, text width=10em, text centered, minimum height=2em]
\tikzstyle{blockns} = [rectangle, draw, text width=5em, text centered, minimum height=5em]
\tikzstyle{blocks} = [rectangle, draw, fill=darkgray!15, text width=5em, text centered, minimum height=5em]
\tikzstyle{block0} = [rectangle, draw, fill=lightgray!15, text width=5em, text centered, minimum height=5em]
\tikzstyle{decision} = [diamond, text width=2cm, minimum width=3cm, minimum height=1cm, text centered, draw, fill=white]
\tikzstyle{line} = [draw, -latex']
\usepackage{makecell}
\usepackage{amssymb}
\usepackage{eurosym}
\usepackage{color}
\usepackage{bm}
\usepackage{booktabs}
\usepackage{multirow}
\usepackage{longtable}
\usepackage{arydshln}
\usepackage{geometry}
\usepackage{enumitem}
\usepackage{algorithm} 
\usepackage{algpseudocode} 
\usepackage{subcaption}
\usepackage{tikz}
\usepackage{pgfplots}
\usetikzlibrary{patterns}
\usetikzlibrary{arrows}
\usepackage{lscape}
\usepackage{amsthm,amsmath,eurosym}
\newtheorem{theorem}{Theorem}
\newtheorem{corollary}{Corollary}
\newtheorem{lemma}{Lemma}

\journal{Computers \& Operations Research}

\usepackage{lineno, blindtext}

\makeatletter
\def\ps@pprintTitle{%
  \let\@oddhead\@empty
  \let\@evenhead\@empty
  \def\@oddfoot{\reset@font\hfil\thepage\hfil}
  \let\@evenfoot\@oddfoot
}
\makeatother

\begin{document}

\begin{frontmatter}
    \title{A Multi-Objective Mixed Integer Linear Programming \\ Model for Thesis Defence Scheduling}
    \author[cegist]{João {\sc  Almeida}\corref{cor1}}\ead{joao.carvalho.almeida@tecnico.ulisboa.pt}
    
    \author[cegist]{Daniel {\sc Rebelo dos Santos}}
    
    \author[cegist]{Jos\'e Rui {\sc  Figueira}}
    \address[cegist]{CEGIST, Instituto Superior T\'{e}cnico,  Universidade de Lisboa, Portugal}

    \cortext[cor1]{Corresponding author at: CEGIST, Instituto Superior T\'{e}cnico,  Universidade de Lisboa, Portugal.}

    \begin{abstract}
    \noindent In this paper, we address the thesis defence scheduling problem, a critical academic scheduling management process, which has been overshadowed in the literature by its counterparts, course timetabling and exam scheduling. Specifically, we address the single defence assignment type of thesis defence scheduling problems, where each committee is assigned to a single defence, scheduled for a specific day, hour and room. We formulate a multi-objective mixed-integer linear programming model, which aims to be a general representation of the problem mentioned above, and that can, therefore, be applied to a broader set of cases than other models present in the literature, which have a focus on the characteristics of their universities. For such a purpose, we introduce a different decision variable, propose constraint formulations that are not regulation and policy specific, and cover and offer new takes on the more common objectives seen in the literature. We also include new objective functions based on our experience with the problem at our university and by applying knowledge from other academic scheduling problems. We also propose a two-stage solution approach. The first stage is employed to find the number of schedulable defences, enabling the optimisation of instances with unschedulable defences. The second stage is an implementation of the augmented $\epsilon-$constraint method, which allows for the search of a set of different and non-dominated solutions while skipping redundant iterations. A novel instance generator for thesis scheduling problems is presented. Its main benefit is the generation of the availability of committee members and rooms in availability and unavailability blocks, resembling their real-world counterparts. A set of 96 randomly generated instances of varying sizes is solved and analysed regarding their relative computational performance, the number of schedulable defences and the distribution of the considered types of iterations. The proposed method can find the optimal number of schedulable defences and present non-dominated solutions within the set time limits for every tested instance.

    \end{abstract}
    \vspace{0.25cm}
\end{frontmatter}

\vfill\newpage

\tableofcontents

\vfill\newpage

\section{Introduction}\label{sec:introduction}
\noindent  The scheduling of thesis defences is a fundamental problem in the academic world. Millions of students worldwide need to prepare and defend their theses in what is often their greatest academic challenge up to that point. Thus, considering the many resource allocation challenges it presents for colleges and universities every year, thesis defence scheduling is one of the most important academic management problems. The literature has established the problem, focusing on each country and university's regulations and culture. Even so, it can be defined through the 6W's question framework: we want to know Who (the committee member) is to be assigned to Which examination committee, performing What role, to Whom (students/thesis defences), When (day and hour) and Where (the specific room). Ultimately, the association of these assignments is called a schedule.

A feasible schedule fulfils a particular set of constraints. For thesis defence scheduling those can be classified under three categories: scheduling complete committees, fulfilling committee composition rules and ensuring committee member and room availability. Moreover, a feasible schedule is not necessarily a ``good'' one. Likewise, two points of view assess the schedules: committee assignment and schedule quality. These points of view are rendered operational by constructing several criteria or objectives to be maximised or minimised. Moreover, regarding committee assignment, the criteria can go from fairly distributing assignments between committee members to assessing their suitability to evaluate the thesis defences. Additionally, the schedule quality objectives might include, for example, promoting compact schedules, satisfying preferred time slot requests from committee members or preventing room changes.

Conversely, many problems may arise when such measures are not considered or when a human scheduler cannot adequately meet them. To begin with, there is the perceived (or actual) bias of the scheduler in favour of some committee members, which might lead to disagreements between the affected parts regarding fairness and transparency considerations. Moreover, if the assigned committee is not knowledgeable enough in the necessary research subjects to assess a thesis, this might lead to inaccurate evaluations of the student's work. Furthermore, as most committee members often have other tasks and are not available at all times, finding a schedule that is as small of an inconvenience for them as possible is not a simple task, which leads to committee members having less available time to dedicate to their research, teaching, or any other activity. Lastly, there is the problem of the scheduling process itself, which just burns too much time for the scheduler, who could put it to more productive use. Likewise, several approaches for dealing with such problems have been considered in the literature. 

Regarding a fair assignment distribution, \cite{KochanikovRudova2013} proposes a constructive heuristic that iteratively assigns defences to committees in such a manner that each committee is appointed the same or at most one more defence than the other committee. Moreover, their improvement phase keeps this distribution intact. \cite{PhamEtAl2015} takes a different approach, minimising the workload differences between the most and least burdened committee members. Conversely, \cite{BattistuttaEtAl2019} is the first to impose an exponential penalty for the number of times a committee member is assigned to a different committee. Finally, \cite{ChristopherWicaksana2021} includes a maximum quota of assignments for each committee member. 

In terms of the methodologies used to assess committee suitability, there is minimal variance. Each consideration of this point of view aims to optimise the matching between expertise or research areas of the committee members and the subject of the defences they are assigned to \citep{HuynhEtAl2012, PhamEtAl2015, BattistuttaEtAl2019, ChristopherWicaksana2021}. Nonetheless, other aspects are also regarded in \cite{ChristopherWicaksana2021}, specifically committee member academic level and previous experience moderating defences. 

As far as reducing the inconvenience caused for committee members, the literature has considered different measures that better suit their situation. \cite{HuynhEtAl2012} considers an objective that minimises the total number of room changes. Moreover, the same work introduces a compactness measure where an increasing number of time slots between two consecutive assignments is penalised. In \cite{BattistuttaEtAl2019}, it is stated that this is a valuable objective. However, due to their structuring of thesis defences into sessions, it is automatically guaranteed. Concerns regarding time slot preferences also arise in \cite{HuynhEtAl2012} and \cite{TawakkalSuyanto2020}, with the latter including an objective specifically regarding students.

There are several types of defence scheduling problems, distinguished by the number of defences that a committee is assigned to: (1) \textit{Single defence assignment} - Each committee is assigned to one defence, which will take place in one slot (day, hour and room). This type of scheduling is addressed in \cite{HuynhEtAl2012, PhamEtAl2015, LimantoEtAl2019, TawakkalSuyanto2020, ChristopherWicaksana2021} and our work; (2) \textit{Session of defences assignment} - Each committee is assigned to a group of defences (session), which will take place in one slot (day, period and room). Evidently, the time slot cannot be just a particular hour in such instances. Instead, it represents an extended period, such as a morning, afternoon, or day. This type of scheduling is addressed in \cite{BattistuttaEtAl2019}; (3) \textit{Hybrid assignment} - A committee is assigned to an extended period, like  in the session of defences assignment type. However, single defences and other required additional committee members are then assigned to an hour within such a period. This type of scheduling is addressed in \cite{KochanikovRudova2013}. 

Several different solution approaches have been applied to the thesis defence scheduling problem, specifically: Mixed-integer linear programming \citep{BattistuttaEtAl2019}; Constraint programming \citep{BattistuttaEtAl2019}; Greedy backtracking hybrid algorithm \citep{SuEtAl2020}; Local search  \citep{KochanikovRudova2013,PhamEtAl2015,BattistuttaEtAl2019, TawakkalSuyanto2020}; Genetic algorithm \citep{HuynhEtAl2012, LimantoEtAl2019}; and Particle swarm optimisation \citep{ChristopherWicaksana2021}. 

In contrast with the most commonly studied academic scheduling problems, specifically exam scheduling and course timetabling \citep{ChaudhuriDe2010}; for a state-of-the-art review of those, we refer the reader to \cite{BabaeiEtAl2015, TanEtAl2021}; and for how inefficient and time-consuming assigning and scheduling committees is for the unfortunate people to whom such tasks are delegated, thesis defence scheduling is remarkably underrepresented in the literature, which still holds several noticeable gaps, with this work being an attempt to fill some of them. In particular, we propose a multi-objective mixed-integer linear programming model, which, to the best of our knowledge, is the first formulation of its type for the single defence assignment type. Moreover, our approach aims to be applicable to a broader range of instances and regulations than the ones found in the literature. Thus, to achieve such a goal, our model includes the following novel characteristics, organised under three levels: 

\begin{itemize}[label={--}]
    \item Main decision variable - our fundamental decision variable, built from the 6W's question framework, is the first to include the committee member's role in each assignment. This is beneficial as it allows for a greater committee composition flexibility than the previous formulations;
    \item Constraints - instead of modelling our university's regulations, we take advantage of the novel decision variable and input the eligible committee members for each role in each thesis defence as a parameter. Therefore, it can then be adapted to fit each instance's needs. Moreover, we are the first to model the problem in a way that considers the possibility of not scheduling every defence, as not all of them might be schedulable, and presenting an ``incomplete'' schedule can still be valuable for the decision-maker;
    \item Objectives - to the best of our knowledge, we are the first to introduce an objective that minimises the number of days a committee member is scheduled to attend a defence. Additionally, we also propose different formulations for several previously defined objectives. Specifically, we introduce a linearisation of an exponential penalty for the number of assignments, compactness, and room change objectives better suited for cases where a committee member is not available for every time slot and a differentiation mechanism for preferred time slots with multiple preference levels. 
\end{itemize}

As for the solution approach, we must be able to find how many thesis defences can be scheduled. Consequently, we propose a two-stage approach, with the first stage being an adaption of our formulation, but with the objective of finding the number of schedulable thesis defences. That number will then be used as a parameter in our model. While new to thesis defence scheduling, similar two-stage approaches have been studied in academic timetabling problems with mixed-integer based solution methods \citep{BurkeEtAl2010Decomposition, SorensenDahms2014, VermuytenEtAl2016}. Their primary motivation was to simplify the search space of an otherwise computationally difficult problem to solve with mixed-integer linear programming. This is not our primary goal, as, instead of using it to solve different objectives separately, we want to find the value of a necessary parameter. Moreover, course timetabling works with heuristic solution methods have also applied two-stage approaches. Similarly to our approach, in \cite{GohEtAl2017, BellioEtAl2021}, the first stage of the heuristic is employed to guarantee the feasibility of the solution. Conversely, other works decompose their heuristics based on different objectives \citep{Santiago-MozosEtAl2005,Al-YakoobSherali2015}. Lastly, while most works regard multiple objectives, none propose approaches that allow the search of the solution space for several solutions, commonly employed in other scheduling problems \citep{LiuEtAl2014, AhmadiEtAl2016, LeiEtAl2016, WangEtAl2017, BorgonjonMaenhout2022}. Contrarily, we adapt the augmented $\epsilon-$constraint method, introduced in \cite{Mavrotas2009} and  \cite{MavrotasFlorios2013}. The main benefits of such an implementation are the assurance that the found solutions  are non-dominated (or Pareto optimal) and the existence of mechanisms to foresee and skip some iterations that will not yield new solutions. 

We propose an instance generator for thesis defence scheduling. It generates the availability of committee members and rooms in availability blocks, resembling their real-world counterparts. Computational experiments are conducted on a set of 96 instances. The main conclusions are that the number of schedulable defences was always relatively easy to optimally determine during the first stage and that, even for the largest considered instances, some non-dominated solutions were found within the set time limits.

The remainder of this paper is organised as follows. Section 2 introduces the multi-objective mixed-integer linear programming model and briefly describes its parameters, variables, constraints, and objective functions. Section 3 is devoted to the approach for identifying non-dominated solutions of the model presented in the previous section. Section 4 addresses the instance generation method. Section 5 presents the computational experiments conducted to test the scalability of the approach and an analysis of the findings. Lastly, Section 6 includes concluding remarks and future research path suggestions.
\section{The Multi-Objective Mixed Integer Linear Programming Model}
\noindent This section presents the multi-objective mixed-integer linear programming model to schedule and assign committees to thesis defences. It introduces the indices of variables and parameters. Then it presents the necessary parameters and variables, followed by a general overview of the objectives and the definition of the constraints. Lastly, the objectives will be revisited and appropriately defined.

\subsection{Indices}\label{SS-Model-Sets}
\noindent The necessary indices for the parameters and variables are presented in this subsection. Let us note that, while they all start with $0$, this value is never used to represent an object. Nonetheless, it is necessary to portray the absence of such objects in specific constraints and objective functions. For example, while there is no committee member $0$, we still need the ability to quantify $0$ committee members.

\begin{itemize}[label={--}]
    \item $i = 0,\ldots, n_i$, are the indices related to the master's thesis defence committee members;
    \item $j = 0,\ldots, n_j$, are the indices related to the master's thesis defences;
    \item $t = 0,\ldots, n_t$, are the indices related to the role of the committee members;
    \item $k = 0,\ldots, n_k$, are the indices related to the days;
    \item $\ell = 0,\ldots, n_\ell$, are the indices related to the available hour slots in each day;
    \item $p = 0,\ldots, n_p$, are the indices related to the available rooms;
    \item $q = 0,\ldots, n_q$, are the indices related to research subjects.  
\end{itemize}

\subsection{Parameters}\label{SS-Model-Parameters}
\noindent The parameters for the model are presented in this subsection. They have been divided into three groups, the first about individual committee members, the second regarding both committee members and rooms, and the last one related to thesis defences.

\begin{enumerate}
    \item \textit{Related to committee members.}
        \begin{itemize}[label={--}]
            \item $e_{ijt} \in \{0,1\}$, is 1 if a committee member $i$ is eligible to be assigned to a certain role $t$ in a designated defence $j$; and 0 otherwise, $i = 1,\ldots, n_i$, for all $j = 1,\ldots, n_j$, $t = 1,\ldots, n_t$;
            \item $c_i\in \mathbb{N}$, is the maximum number of committees a committee member $i$ can be assigned to, for all $i = 1,\ldots, n_i$;
            \item $u_i \in \mathbb{N}$, is the weight assigned to each committee member $i$, for all $i = 1,\ldots, n_i$;
            \item $l_{ik\ell} \in \mathbb{N}_0$, is the preference level that each committee member $i$ holds for every time slot (day $k$ and hour $\ell$), 0 represents unavailability, for all $i = 1,\ldots, n_i$, $k = 1,\ldots, n_k$, $\ell = 1,\ldots, n_\ell$; 
            \item $r_{iq} \in \{0,1\}$, is 1 if a committee member $i$ has knowledge regarding the research subject $q$; and 0 otherwise, for all $i = 1,\ldots, n_i$, $q = 1,\ldots, n_q$;
        \end{itemize}
    \item \textit{Related to committee members and rooms.}
        \begin{itemize}[label={--}]
            \item $b_i \in \{0,\ldots, d - 1\}$, is the number of hour slots following the end of a defence in which a committee member $i$ would consider the scheduling of a different one as compact, in our formulation we considered it as always being smaller than the duration  $d$ of a defence, for all $i = 1,\ldots, n_i$;
            \item $v_{i\ell} \in \{0, \dots, n_{v_i}\}$, is the weight given by a committee member $i$ to each hour slot $\ell$ considered compact, for all $i = 1,\ldots, n_i$, $\ell = 0,\dots b_i$;
            \item $a_i \in \{0,\ldots, d - 1\}$, is the number of hour slots following the end of a defence in which a committee member $i$ would consider changing rooms problematic, for all $i = 1,\ldots, n_i$;
            \item $h_{i\ell} \in \{0, \dots, n_{h_i}\}$, is the penalty given by a committee member $i$ after a room change to each hour slot $\ell$ considered problematic, for all $i = 1,\ldots, n_i$, $\ell = 0,\dots, a_i$; 
            \item $m_{k\ell p} \in \{0,1\}$, is 1 if a room $p$ is available at a certain time slot ($k$, $\ell$); and 0 otherwise, for all $k = 1, \ldots, n_k$, $\ell = 1, \ldots, n_\ell$, $p = 1, \ldots, n_p$;
        \end{itemize}
    \item \textit{Related to thesis defences.}
        \begin{itemize}[label={--}]
            \item $d \in \mathbb{N}$, is the length (duration) of a defence in hour slots;
           
            \item $\overline{t}_{jq}\in \{0,1\}$, is 1 if a defence $j$ studies a certain research subject $q$; and 0 otherwise, for all $j = 1,\ldots, n_j$, $q = 1,\ldots, n_q$;
            \item $g \in \{0, \ldots, n_j\}$, is the number of complete committees (thesis defences) to be scheduled. If one knows that all defences are schedulable, then $g = n_j$. Otherwise, finding $g$ becomes part of the problem itself.
        \end{itemize}
\end{enumerate}

\subsection{Variables}\label{SS-Model-Variables}
\noindent The definition of the decision and auxiliary variables is presented in this subsection. The first is built from the 6W's question framework and represents a committee member's assignment to a role within a thesis defence, happening at a specific time and room. The others are used to define concepts such as the number of scheduled days, workloads, defence research subject coverage by its committee, compactness and room change penalties, and assigned roles within a committee.

\begin{enumerate}
    \item \textit{Decision variables.}
        \begin{itemize}[label={--}]
            \item $x_{ijtk\ell p}\in \{0,1\}$, which is equal to $1$ if a committee member $i$ is assigned to thesis defence $j$, performing a role $t$, in day $k$, at hour slot $\ell$ and in room $p$; and $0$ otherwise, for all $i = 1,\ldots, n_i$, $j = 1,\ldots, n_j$, $t = 1,\ldots, n_t$, $k = 1,\ldots, n_k$, $\ell = 1,\ldots, n_\ell$, and $p = 1,\ldots, n_p$; 
        \end{itemize}
     \item \textit{Auxiliary variables.}
        \begin{itemize}[label={--}]  
            \item $y_{jk\ell p}\in \{0,1\}$, is 1 if the thesis defence $j$ is scheduled for day $k$, at hour $\ell$ in room $p$; and $0$ otherwise, for all $j = 1,\ldots, n_j$, $k = 1,\ldots, n_k$, $\ell = 1,\ldots, n_\ell$, and $p = 1,\ldots, n_p$; 
            \item $\overline{y}_{ik\ell p}\in \{0,1\} $, is $1$ if committee member $i$ is assigned to any defence on day $k$, hour $\ell$ and room $p$; and $0$ otherwise, for all $i = 1,\ldots, n_i$, $k = 1,\ldots, n_k$, $\ell = 1,\ldots, n_\ell$, and $p = 1,\ldots, n_p$; 
            \item$\hat{y}_{ijk}\in \{0,1\} $, is $1$ if committee member $i$ is assigned to $j$ defences on day $k$; and $0$ otherwise, for all $i = 1,\ldots, n_i$, $j = 0, \ldots, n_j$ and $k = 1,\ldots, n_k$; 
            \item $w_{ij}\in \{0,1\} $, is $1$ if committee member $i$ is assigned to $j$ thesis defences; and $0$ otherwise, for all $i = 1,\ldots, n_i$, and $j = 0,\ldots, n_j$; 
            \item $\overline{w}_{ik}\in \{0,1\}$, is $1$ if committee member $i$ is assigned to thesis defences in $k$ days; and $0$ otherwise, for all $i = 1,\ldots, n_i$, and $k = 0,\ldots, n_k$; 
            \item $s_{ijq} \in \{0, 1\}$, is 1 if $i$ committee members of the committee of defence $j$ have $q$ research subject in common with said defence; and 0 otherwise, for all $i = 0, \ldots, n_i$, $j = 1,\ldots, n_j$, and $q = 1,\ldots, n_q$;
            \item $\overline{s}_{ik\ell p} \in \mathbb{N}_{0}$, is the compactness value for committee member $i$, scheduled to attend any defence at day $k$, hour $\ell$ and room $p$, for all  $i = 1,\ldots, n_i$, $k = 1,\ldots, n_k$, $\ell = 1,\ldots, n_\ell$, and $p = 1,\ldots, n_p$;
            \item  $\hat{s}_{ik\ell p} \in \mathbb{N}_{0}$, is a room change penalty for committee member $i$, scheduled to attend any defence at day $k$, hour $\ell$ and room $p$, for all  $i = 1,\ldots, n_i$, $k = 1,\ldots, n_k$, $\ell = 1,\ldots, n_\ell$, and $p = 1,\ldots, n_p$;
        \end{itemize}
\end{enumerate}

\subsection{Objective functions}\label{SS-Model-Objectives}
\noindent The objective functions are presented in this subsection. Two points of view are used to assess the assignment quality. The first one, committee assignment, is related to both fair workload distribution and the matching between the expertise of the committee members and the defences they are assigned to. The second one, schedule quality, is related to the assessment of the quality of the schedules of individual committee members. Moreover, the objective functions used to render them operational in assessing the quality of each feasible schedule/solution are introduced here. Nonetheless, further constraints and variables need to be defined before the mathematical expressions for each objective can be specified. In Subsection \ref{SS-Model-BackObjectives} we will return to the objective functions and properly define them. Furthermore, the additional constraints are listed in Subsection \ref{SS-Model-Constraints}.

    \begin{enumerate}
        \item \textit{Point of view of committee assignment.} This point of view is operationalised by the following objectives. 
            \begin{enumerate}
                \item \textit{Minimise workloads}. The workload is the number of committee assignments a committee member has. To ensure fairness, this number is squared. We want to achieve a balanced workload distribution between committee members with this criterion.
                \item \textit{Maximise research subject coverage}. Research subject coverage is the percentage of research subjects in a defence that its committee covers. We want to maximise such a coverage with this criterion.  
                \item \textit{Maximise committee member suitability}. Committee member suitability is defined as the number of research subjects each committee member has in common with their assigned defences. We want to maximise such suitability with this criterion.
            \end{enumerate}
        \item \textit{Point of view of schedule quality.} This point of view is operationalised by the following objectives. 
            \begin{enumerate}
                \item \textit{Minimise non-consecutive assignments}. Each assignment is given a compactness value. This value is the product of the committee member's weight and the weight for the interval in which a defence is scheduled. The committee member defines the latter regarding their preferences over time intervals between defences. We want to minimise the difference between the highest potential compactness value for an assignment and its effective value, producing more compact schedules according to the committee members' preferences with this criterion.
            \item \textit{Minimise the non-satisfaction of time slot preferences}. A penalty is given whenever a committee member is assigned one of their stated undesirable time slots. We want to minimise the occurrence of such assignments with this criterion. 
            \item \textit{Minimise committee days}. Committee days are the number of days a committee member is scheduled to attend a defence. To ensure fairness, this number is squared. We want to minimise such a number with this criterion.
            \item \textit{Minimise room changes}. Each committee member defines a time frame considered problematic for a room change between defences, penalising such events. We want to promote room stability with this criterion.
        \end{enumerate}
    \end{enumerate}

\subsection{Constraints}\label{SS-Model-Constraints}
\noindent The necessary constraints to model the feasible region of the problem are presented in this subsection. They fall under four categories. The first concerns the scheduling of complete committees or thesis defences. The second regards the respect for committee member assignment rules. The third guarantees that committee members and rooms are available for their corresponding assignments. Finally, the fourth defines the values for several auxiliary variables present in the objective functions.

\begin{enumerate}
    \item \textit{Scheduling complete committees.} These constraints define a complete committee and ensure that
every schedulable defence is assigned one.
        \begin{enumerate}
            \item \textit{Complete committee definition}. A complete committee is a set of $n_t$ assignments for a defence, $j$, all with a different appointed role, $t$, in the same slot, (day $k$, hour $\ell$, and room $p$). Moreover, for a defence to occur, it must have such a committee assigned to it. Conversely, no assignment that is not incorporated into one can exist, as it would occupy a slot that is not being used. Likewise, this constraint defines the auxiliary variable $y_{jk\ell p}$, which takes the value $1$ if a defence, $j$, has a complete committee assigned to it and $0$ otherwise. Thus, the left-hand-side of the equation can also only take binary values. The sum on the mentioned side represents the number of committee members, $i$, assigned to a defence, $j$, to perform a role, $t$, on a given slot, ($k, \ell, p$). Consequently, since this sum can be at most $1$, each role, $t$, can only be filled once in a complete committee. Moreover, evidently, $y_{jk\ell p}$ can only take one value for a given defence, $j$, and slot, ($k, \ell, p$), ergo, the left-hand-side of the equality can also only take one value for the same defence, $j$, slot, ($k, \ell, p$) and for every role, $t = 1, \ldots, n_t$, meaning that for $ t \neq \overline{t}$, the following equality  must be verified: $\sum_{i = 1}^{n_i} x_{ijtk\ell p} = \sum_{i = 1}^{n_i} x_{ij\overline{t}k\ell p}$. Accordingly, if a role, $t$, is assigned to a defence, $j$, any other role, $\overline{t}$, must also be assigned. Interchangeably, if any role, $t$, is not assigned, any other role, $\overline{t}$, must also not be assigned. Therefore, either a defence is assigned to a complete committee or no assignments for such a defence can occur.

                   \begin{equation}\label{Const - comm size - eq1}
                    {\displaystyle \sum_{i = 1}^{n_i} x_{ijtk\ell p} =  y_{jk\ell p}, \;\;\; j = 1, \ldots, n_j,\; t = 1, \ldots, n_t \; k = 1, \ldots, n_k,\; \ell = 1, \ldots, n_\ell, \; p = 1, \ldots, n_p}
                \end{equation}

            \item \textit{Single committee assignment}. If a defence, $j$, can be scheduled, it should only be assigned one committee and appointed one slot, (day $h$, hour $\ell$, and room $p$). Thus, in this constraint, we state that for a defence, $j$, the number of complete committees assigned to it is less or equal to $1$.
  
                \begin{equation}\label{const - defences scheduled once - eq 4}
                    \displaystyle \sum_{k = 1}^{n_k}\sum_{\ell = 1}^{n_\ell}\sum_{p = 1}^{n_p} y_{jk\ell p} \leqslant 1, \;\;\;  j = 1, \ldots, n_j
                \end{equation}
            \item \textit{Complete committees (thesis defences) to be scheduled}. In each instance of the thesis defence scheduling problem, a defined number of committees (thesis defences) can be assigned and scheduled, denoted by $g$. If one knows that all the defences can be scheduled, then this number is  the number of defences, i.e., $g = n_j$. However, some defences may not be schedulable due to conflicting committee member availabilities, lack of enough eligible committee members, lack of rooms, or others. In such cases, finding the value for $g$ becomes an indispensable part of the problem. In this constraint, assuming the value of $g$ is already known, we enforce the number of assigned complete committees (thesis defences) as the number of schedulable complete committees, $g$.
                \begin{equation}\label{const - defences scheduled once - eq 5}
                    {\displaystyle \sum_{j = 1}^{n_j}\sum_{k = 1}^{n_k}\sum_{\ell = 1}^{n_\ell}\sum_{p = 1}^{n_p} y_{jk\ell p} = g}
                \end{equation}

        \end{enumerate}
    \item \textit{Committee Composition.} These constraints ensure the eligibility of the committee members to
perform their assignments.
        \begin{enumerate}

            \item \textit{Committee member eligibility}. Different universities and their departments have distinct regulations for the eligibility of committee members, $i$, to perform specific roles, $t$, within each committee for a defence, $j$. Thus, we do not attempt to include such rules within our model. Conversely, we aggregate them in a parameter, $e_{ijt}$, which takes the value 1 if a committee member, $i$, is eligible to perform a role, $t$, in a defence, $j$, and 0 otherwise. In this constraint, we state that if a given committee member, $i$, is non-eligible to perform a role, $t$, in the committee of a defence, $j$, that is, if $e_{ijt} = 0$, then no assignment that involves such a combination can occur, that is, the left-hand-side of the equation must also be 0. Contrarily, if such a combination is possible, that is, if $e_{ijt} = 1$, then the equality still holds, as, logically, a committee member, $i$, can be assigned at most once to a given defence, $j$.

                \begin{equation} \label{const - eligible members}
                    \displaystyle \sum_{k = 1}^{n_k}\sum_{\ell = 1}^{n_\ell}\sum_{p = 1}^{n_p} x_{ijtk\ell p} \leqslant e_{ijt}, \;\;\;  i = 1,\ldots, n_i,\; j = 1,\ldots, n_j,\; t = 1,\ldots, n_t
               \end{equation}
            \item \textit{Maximum number of committees assigned to a committee member}. This committee member eligibility requirement cannot be represented by the eligibility parameter, $e_{ijt}$. Thus, we included it as another constraint. In cases where there is no such regulation, the value for the maximum number of committees assigned to a committee member, $i$, represented by $c_i$, is equal to the number of thesis defences, that is, $c_i = n_j$. In this constraint, we ensure that the sum of the assignments, which occur when $x_{ijtk\ell p} = 1$, for a committee member, $i$, does not exceed the maximum allowed number of committees, $c_i$, for that committee member.
                \begin{equation}\label{const - max committees}
                     \displaystyle \sum_{j = 1}^{n_j}\sum_{t = 1}^{n_t}\sum_{k = 1}^{n_k}\sum_{\ell = 1}^{n_\ell}\sum_{p = 1}^{n_p} x_{ijtk\ell p} \leqslant  c_i, \;\;\;  i = 1, \ldots, n_i
                \end{equation}
        \end{enumerate}
    \item \textit{Committee member and room availability}. These constraints guarantee that committee members
and rooms are available for each assignment.
        \begin{enumerate}
            \item \textit{Committee member time slot availability}. Committee members have different obligations other than attending thesis defences. Consequently, they are not available to be assigned to every time slot. The committee member availability parameter, $l_{ik \ell}$, takes a value greater than or equal to $1$ if a committee member, $i$, is available to be assigned at a day, $k$, and an hour, $\ell$, and 0 if they are not. In this, we state that if a given committee member, $i$, is not available at a day, $k$, and an hour, $\ell$, that is, if  $l_{ik\ell} = 0$, then no assignment that involves such a combination can occur, that is, the left-hand-side of the equation must also be 0. Contrarily, if such a combination is possible, that is, if $l_{ik\ell} \geqslant 1$, then the equality still holds, as, logically, a committee member, $i$, can only be given at most one assignment at a particular time slot, ($k, \ell$).
                \begin{equation}\label{const - member availability}
                    \displaystyle \sum_{j = 1}^{n_j}\sum_{t = 1}^{n_t}\sum_{p = 1}^{n_p} x_{ijtk\ell p} \leqslant  l_{ik\ell}, \;\;\; i = 1,\ldots, n_i,\;
                    k = 1,\ldots, n_k,\; \ell = 1,\ldots, n_\ell
                \end{equation}
            \item \textit{Committee member assignment juxtaposition}. A committee member, $i$, cannot be assigned to more than one defence, $j$, starting at a day, $k$, and an hour, $\ell$. Moreover, that committee member is also unavailable to attend any other defence that begins at any given point before the end of such a defence, $j$. In other words, until the hour $\ell + d$ is reached on the same day, $k$. Thus, in this constraint, we ensure that if there is an assignment, $x_{ijtk\overline{\ell} p} = 1$, for a committee member, $i$, in a day, $k$, at an hour, $\overline{\ell}$, there cannot be any other assignment for the same committee member, in an hour that occurs before the end of the previous defence, that is, in any hour between and including $\overline{\ell}$ and $\overline{\ell} + d - 1$. Additionally, this constraint also ensures that a committee member is not assigned more than one role, $t$, in a defence, $j$, as that would mean that said committee member would have two different assignments in the same time slot, ($k, \ell$).
                \begin{equation}\label{const - member juxtaposition}
                    \displaystyle \sum_{j = 1}^{n_j}\sum_{t = 1}^{n_t}\sum_{\ell = \overline{\ell}}^{\overline{\ell} + d - 1}\sum_{p = 1}^{n_p} x_{ijtk\ell p} \leqslant 1, \;\;\; i = 1, \ldots, n_i,\; k = 1, \ldots, n_k,\; \overline{\ell} = 1, \ldots, n_\ell - d + 1
                \end{equation}

            \item \textit{Room time slot availability}. A room's purpose might not just be hosting thesis defences. Thus, it is natural that it happens to be booked for any other event at some point. The room availability parameter, $m_{k\ell p}$, takes the value $1$ if a room, $p$, is available to host a defence at a day, $k$, and an hour, $\ell$, and $0$ otherwise. In this constraint, we state that, for a given slot, ($k, \ell, p$), the sum of its assigned complete committees, $y_{jk\ell p} = 1$,  always takes a value lower or equal to that of $m_{k \ell p}$. Accordingly, whenever a room is unavailable, that is, $m_{k \ell p} = 0$, it cannot host any defence, and this sum is correctly set to $0$. Moreover, if the room is available, that is, $m_{k \ell p} = 1$, then, at most, one defence, $j$, can be assigned such a slot, ($k, \ell, p$).
                \begin{equation}\label{const - room availability}
                    \displaystyle \sum_{j = 1}^{n_j}y_{jk\ell p} \leqslant m_{k\ell p},\;\;\;
                    k = 1,\ldots, n_k,\; \ell = 1,\ldots, n_\ell,\; p = 1,\ldots, n_p
                \end{equation}
            \item \textit{Room capacity}. In our formulation, we considered that a room could not hold more than one defence at a time. Variable $y_{jk\ell p}$ takes the value 1 if a defence, $j$, is assigned to a day, $k$, an hour, $\ell$, and a room, $p$. Consequently, for the same day, room and between the start of defence $j$, at hour $\ell$, and its end, at hour $\ell + d$, the point at which the room can be scheduled again, there cannot be more than one $y_{jk\ell p} = 1$. In this constraint, this is achieved by stating that, for any day, $k$, hours between, and including, $\ell$ to $\ell + d - 1$, and a room, $p$, the sum of the values of the variable $y_{jk\ell p}$ must be less or equal to $1$.
                \begin{equation}\label{const - room capacity}
                    \displaystyle \sum_{j = 1}^{n_j}\sum_{\ell = \overline{\ell}}^{\overline{\ell} + d - 1} y_{jk\ell p} \leqslant 1, \;\;\; k = 1, \ldots, n_k,\; \overline{\ell} = 1, \ldots, n_\ell - d + 1,\; p = 1, \ldots, n_p
                \end{equation}
        \end{enumerate}
        \item \textit{Objective functions measures.} These constraints define the values for the auxiliary variables necessary for some of the objective functions.
            \begin{enumerate}
                \item \textit{Research subject coverage definition.} For a defence, $j$, research subject coverage is the percentage of its studied research subject covered by the areas of expertise of its committee members. To implement such an objective, we first need to define an auxiliary binary variable, $s_{ijq}$, which takes the value 1 if a defence, $j$, which studies a research subject, $q$, that is, $\overline{t}_{jq} = 1$, has a number, $i$, of committee members assigned to its committee, who have said subject as one of their areas of expertise. Let us note that, whenever a committee member, $i$, studies a research subject, $q$, then $r_{iq} = 1$. The value for such a variable is defined in Constraint \eqref{const - res perc - eq 1} by stating that the product of $s_{ijq}$ by the number, $i$, of committee members assigned to that defence is equal to the sum of the assignments, $x_{ijtk\ell p} = 1$, where the research subject, $q$, is in both the studied subjects of the defence and the areas of expertise of the committee member, that is $r_{iq}\overline{t}_{jq} = 1$. Furthermore, a defence, $j$, cannot be assigned to more than one number of committee members with a research subject, $q$, in common with it. Thus, with Constraint \eqref{const - res perc - eq 2}, we ensure that this value is unique for each combination of defence, $j$, and research subject, $q$.

                \begin{equation}\label{const - res perc - eq 1}
                    \displaystyle\sum_{i = 0}^{n_i} is_{ijq} =  \sum_{i = 1}^{n_i}\sum_{t = 1}^{n_t}\sum_{\ell = 1}^{n_\ell}\sum_{p = 1}^{n_p}  r_{iq}\overline{t}_{jq}x_{ijtk\ell p}, \;\;\; j = 1, \ldots, n_j, \; q = 1, \ldots, n_q
                \end{equation}
                \begin{equation}\label{const - res perc - eq 2}
                    \displaystyle \sum_{i = 0}^{n_i} {s}_{ijq} = 1, \;\;\; j = 1,\ldots, n_j, \; q = 1, \ldots, n_q
                \end{equation}

                \item \textit{Compactness value definition.} We defined a compact assignment of a committee member, $i$, to a day, $k$, at an hour, $\ell$, as one that occurs within a specific time frame, $b_i$, after the end of a different assignment for such a committee member. That is, if a committee member, $i$, is assigned to a defence, $j$, at a day, $k$, and an hour, $\ell$, this assignment is considered compact on the condition that the same committee member is assigned to a different defence, $\overline{j}$, in the same day, $k$, between hours $\ell - d$ and $\ell - d - b_i$. Moreover, the parameter $v_{i\ell}$, distinguishes the hour slots within such a time frame, as they might have different perceived values for a committee member, $i$. Thus, the compactness value for an assignment, $\overline{s}_{ik \ell}$, is 0 if a committee member, $i$, does not have a different assignment ending between hours $\ell$ and $\ell - b_i$, or $v_{i \overline{\ell}}$ if he does have such an assignment ending at $\overline{\ell}$ hour slots before the start of the new assignment a different hour slot, $\ell$.
                
                \hspace{0.5cm} Before we define the compactness variable, $\overline{s}_{ik \ell}$, we define a different variable, $\overline{y}_{ik\ell p}$, which takes the value $1$ if a committee member, $i$, is assigned to any defence at a day, $k$, an hour, $\ell$ and a room, $p$, and $0$ otherwise. This is done in Constraint \eqref{const - comp score - eq 2}. Moreover, we denote that the right-hand-side of the equality in the constraint above never takes a value greater than $1$, as a committee member, $i$, cannot be assigned more than one defence or role in the same slot, ($k, \ell, p$).
                
                \hspace{0.5cm} With this new variable, we can define the compactness variable, $\overline{s}_{ik\ell}$, as the product between the sums, $\sum_{p = 1}^{n_p} \overline{y}_{ik\ell p}$, and, $\sum_{\overline{\ell} = 0}^{b_i} \sum_{p = 1}^{n_p} v_{i\overline{\ell}}\overline{y}_{ik\hat{\ell}p}$, with $\hat{\ell} = \ell - d - \overline{\ell}$. Nonetheless, this would impair the linearity of the model. Thus, to linearise the aforementioned product, we opted for a big-$M$ formulation, with $M = n_{vi}$, which is the highest value the compactness variable, $\overline{s}_{ik\ell}$, can take for a committee member, $i$. This formulation is represented in the following constraints.

                \begin{equation}\label{const - comp score - eq 2}
                    \displaystyle \overline{y}_{ik\ell p} = \sum_{j = 1}^{n_j}\sum_{t = 1}^{n_t}x_{ijtk\ell p}, \;\;\; i = 1,\ldots, n_i,\;
                    k = 1,\ldots, n_k,\; \ell = 1,\ldots, n_\ell,\; p = 1,\ldots, n_p
                \end{equation}
 
                \begin{equation}\label{const - comp score - eq 4}
                    \displaystyle \overline{s}_{ik\ell} \geqslant 0, \;\;\; i = 1,\ldots, n_i,\;
                    k = 1,\ldots, n_k,\; \ell = 1,\ldots, n_\ell
                \end{equation}
                \begin{equation}\label{const - comp score - eq 5}
                    \displaystyle \overline{s}_{ik\ell} \leqslant n_{v_i}\sum_{p = 1}^{n_p} \overline{y}_{ik\ell p}\;\;\; i = 1,\ldots, n_i,\;
                    k = 1,\ldots, n_k,\; \ell = 1,\ldots, n_\ell
                \end{equation}
                \begin{equation}\label{const - comp score - eq 6}
                    \displaystyle \overline{s}_{ik\ell} \leqslant \sum_{\overline{\ell} = 0}^{b_i} \sum_{p = 1}^{n_p} v_{i\overline{\ell}}\overline{y}_{ik\hat{\ell}p},\;\;\; i = 1,\ldots, n_i,\;
                    k = 1,\ldots, n_k,\; \ell = d,\ldots, n_\ell, \; \hat{\ell} = \ell - d - \overline{\ell}
                \end{equation}
                \begin{equation}\label{const - comp score - eq 7}
                    \begin{array}{ll}
                        \displaystyle \overline{s}_{ik\ell} \geqslant \sum_{\overline{\ell} = 0}^{b_i} \sum_{p = 1}^{n_p} v_{i\overline{\ell}}\overline{y}_{ik\hat{\ell}p} - n_{v_i}\left(1 - \sum_{p = 1}^{n_p} \overline{y}_{ik\ell p}\right), & i = 1,\ldots, n_i,\;
                    k = 1,\ldots, n_k, \\& \ell = d,\ldots, n_\ell, \;\hat{\ell} = \ell - d - \overline{\ell} \\
                    \end{array}
                \end{equation}
                \item \textit{Workload definition.} The workload for a committee member, $i$, is defined as the number, $j$, of committees they are assigned to. It would be possible to represent it as an integer variable. Still, in such a case, it would not be possible to consider its square in the objective function while keeping its linearity. However, the exponential penalty in the objective function can be linearised by representing it through a variable, $w_{ij}$, which takes the value 1 if a committee member, $i$, is assigned to a number, $j$, of committees, and 0 otherwise. The value for such a variable is defined in Constraint \eqref{const - wl def - eq1} by stating that the product of $w_{ij}$ by the number, $j$, of defences assigned to a committee member, $i$, is equal to the sum of the assignments, $x_{ijtk\ell p} = 1$, for that same committee member. Furthermore, a committee member, $i$, cannot be assigned to more than one number of defences, $j$. Thus, with Constraint \eqref{const - wl def - eq2}, we ensure that this value is unique for each one.  
                    \begin{equation}\label{const - wl def - eq1}
                        \displaystyle \sum_{j = 0}^{c_i} jw_{ij} =  \sum_{j = 1}^{n_j}\sum_{t = 1}^{n_t}\sum_{k = 1}^{n_k}\sum_{\ell = 1}^{n_\ell}\sum_{p = 1}^{n_p} x_{ijtk\ell p}, \;\;\; i = 1, \ldots, n_i
                    \end{equation}
                    \begin{equation}\label{const - wl def - eq2}
                        \displaystyle \sum_{j = 0}^{c_i} w_{ij} = 1, \;\;\; i = 1, \ldots, n_i
                    \end{equation}
                \item \textit{Committee days definition.} A committee day is defined as a day when a committee member has a defence scheduled. To represent this concept, we introduce a variable, $\hat{y}_{ijk}$, which takes value 1 if a committee member, $i$, is assigned to a number, $j$, of committees in a day, $k$, and 0 otherwise. To define such a variable, in Constraint \eqref{const - comm days - eq 1}, we define its value and in \eqref{const - comm days - eq 2} we ensure its uniqueness for each combination of committee member, $i$, and day, $k$, similarly to the definitions of variables $s_{ijq}$, in Constraints \eqref{const - res perc - eq 1} and \eqref{const - res perc - eq 2}, and $w_{ij}$, in Constraints \eqref{const - wl def - eq1} and \eqref{const - wl def - eq2}. To institute fairness in the distribution of committee days, an exponential penalty on the number for each committee member is included in their respective objective function. However, to keep the model linear, we still need another variable, $\overline{w}_{ik}$, which takes the value $1$ if a committee member, $i$, has defences scheduled on a number, $k$, of days, and 0 otherwise. The value for such a variable is defined in Constraint \eqref{const - comm days - eq 4}, and, in Constraint \eqref{const - comm days - eq 5}, we ensure that this value is unique for each committee member. Moreover, let us note that, while similar, the latter two constraints have a noteworthy difference when compared to the other constraints referenced in this point. Specifically, the right-hand-side of Constraint \eqref{const - comm days - eq 4}, does not include a variation of the sum of the assignments, $x_{ijtk\ell p} = 1$, involving instead another variable, $\hat{y}_{ijk}$, moreover, while this variable, $\hat{y}_{ijk}$, is defined for a $j = 0, \ldots, n_j$, the sum must start in $j = 1$, as we do not want to count the days, $k$, where a committee member, $i$, has 0 defences assigned.

                \begin{equation}\label{const - comm days - eq 1}
                    \displaystyle\sum_{j = 0}^{c_i} j\hat{y}_{ijk} =  \sum_{j = 1}^{n_j}\sum_{t = 1}^{n_t}\sum_{\ell = 1}^{n_\ell}\sum_{p = 1}^{n_p} x_{ijtk\ell p}, \;\;\; i = 1, \ldots, n_i, \; k = 1, \ldots, n_k
                \end{equation}
                \begin{equation}\label{const - comm days - eq 2}
                    \displaystyle \sum_{j = 0}^{n_j} \hat{y}_{ijk} = 1, \;\;\; i = 1,\ldots, n_i, \; k = 1, \ldots, n_k
                \end{equation}
                \begin{equation}\label{const - comm days - eq 4}
                    \displaystyle \sum_{k = 0}^{n_k} k\overline{w}_{ik} = \sum_{j = 1}^{c_i}\sum_{k = 1}^{n_k} \hat{y}_{ijk} , \;\;\; i = 1, \ldots, n_i
                \end{equation}

                \begin{equation}\label{const - comm days - eq 5}
                    \displaystyle \sum_{k = 0}^{n_k} \overline{w}_{ik} = 1, \;\;\; i = 1, \ldots, n_i
                \end{equation} 

                \item \textit{Room change penalty definition.} A room change is considered problematic if a committee member, $i$, is not given a certain amount of time, $a_i$, between the end of an assignment, $\overline{j}$, and the beginning of another, $j$, which is scheduled for a different room, $\overline{p}$, than the first one, $p$. Moreover, parameter $h_{i\ell}$ distinguishes the hour slots within such a time-frame, as they might have different perceived penalties for a committee member, $i$. The room change variable, $\hat{s}_{ik\ell p}$, is the variable that represents the room change penalty that the assignment of a committee member, $i$, to a day, $k$, an hour, $\ell$ and a room, $p$, would incur. This variable can be defined as the product between the variable $\overline{y}_{ik\ell p}$, which takes the value 1 if a committee member, $i$, is assigned to any defence at a day, $k$, an hour, $\ell$, and a room, $p$, and the sum $\sum_{\overline{\ell} = 0}^{a_i} \sum_{\overline{p} = 1}^{n_p} h_{i\overline{\ell}}\overline{y}_{i k \hat{\ell}\overline{p}}$, with $\hat{\ell} = \ell - d - \overline{\ell}$, which will take the value of parameter $h_{i\overline{\ell}}$ if a committee member, $i$, is assigned to a different defence, $\overline{j}$, in the same day, $k$, in hour $\ell - d - \overline{\ell}$, which is allocated a different room, $\overline{p}$. However, this product would not be linear. Thus, we opted for a big-$M$ formulation, with the big-$M$ being bounded by the highest value parameter $h_{i\ell}$ can take, that is, $M = n_{h_i}$ to linearise said product. This formulation is represented in the following constraints.
                \begin{equation}\label{const - room - eq 1}
                    \displaystyle \hat{s}_{ik\ell p} \geqslant 0, \;\;\; i = 1,\ldots, n_i,\;
                    k = 1,\ldots, n_k,\; \ell = 1,\ldots, n_\ell,\; p = 1, \ldots, n_p
                \end{equation}
                \begin{equation}\label{const - room - eq 2}
                    \displaystyle \hat{s}_{ik\ell p} \leqslant n_{h_i} \overline{y}_{ik\ell p}\;\;\; i = 1,\ldots, n_i,\;
                    k = 1,\ldots, n_k,\; \ell = 1,\ldots, n_\ell,\; p = 1, \ldots, n_p
                \end{equation}
                \begin{equation}\label{const - room - eq 3}
                \begin{array}{ll}

                    \displaystyle \hat{s}_{ik\ell p} \leqslant \sum_{\overline{\ell} = 0}^{a_i} \sum_{\overline{p} = 1}^{n_p} h_{i\overline{\ell}}\overline{y}_{i k \hat{\ell} \overline{p}},& i = 1,\ldots, n_i,\;
                    k = 1,\ldots, n_k, \;\ell = d,\ldots, n_\ell,\\& \hat{\ell} = \ell - d - \overline{\ell}, \; p = 1, \ldots, n_p, \;p \neq \overline{p}\\
                \end{array}
                \end{equation}
                \begin{equation}\label{const - room - eq 4}
                    \begin{array}{ll}
                        \displaystyle \hat{s}_{ik\ell p} \leqslant \sum_{\overline{\ell} = 0}^{a_i} \sum_{\overline{p} = 1}^{n_p} h_{i\overline{\ell}}\overline{y}_{i k \hat{\ell}\overline{p}} - n_{h_i}(1 -  \overline{y}_{ik\ell p}), & i = 1,\ldots, n_i,\;
                    k = 1,\ldots, n_k, \\& \ell = d,\ldots, n_\ell, \; \hat{\ell} = \ell - d - \overline{\ell},\\
                    & p = 1, \ldots, n_p,\; p \neq \overline{p}
                    \end{array}
                \end{equation}
            \end{enumerate}
    \end{enumerate}

\subsection{Back to the objective functions}\label{SS-Model-BackObjectives}
\noindent In this subsection, now that we have defined all the necessary constraints, the objective functions within the two points of view are revisited and adequately determined.
    \begin{enumerate}
        \item \textit{Point of view of committee assignment}. This point of view is operationalised by the following objectives.
            \begin{enumerate}
                \item \textit{Minimise workloads}. This criterion considers a product, $u_{i}j^{2}w_{ij}$. The variable, $w_{ij}$, takes the value 1 if a committee member, $i$, is assigned to a number, $j$, of defences, and 0 otherwise. By multiplying it by $j^2$, we can keep the linearity of the model while applying an exponential penalty representing the square of the number of defences a committee member, $i$, is assigned to, promoting fairness in workload distribution. Moreover, the committee member's weight, $u_i$, is also considered. We want to promote a fair distribution of workloads between committee members by minimising this objective.
                    \begin{equation}\label{ob1}
                        {\displaystyle \min z_1(w) = \sum_{i = 1}^{n_i}\sum_{j = 1}^{n_j}u_{i}j^{2}w_{ij}}
                    \end{equation}
                \item \textit{Maximise research subject coverage}. The research subject coverage for a defence, $j$, is computed through a quotient between the number of its research subjects, $q$, covered by its committee as the numerator and the sum of all of its research subjects as the numerator. Let us note that, similarly to Constraint \eqref{const - comm days - eq 4}, while variable $s_{ijq}$ is defined for $i = 0, \ldots, n_i$, the sum must start on $i = 1$, as we do not want to include the research subjects that are covered 0 times by the defence's assigned committee. We want to maximise the sum of the coverages of all defences with this objective.
                    \begin{equation}
                        {\displaystyle \max z_2(s) = \left(\sum_{j = 1}^{n_j}\sum_{q = 1}^{n_q}\overline{t}_{jq}\right)^{-1}  \sum_{i = 1}^{n_i} \sum_{j = 1}^{n_j}\sum_{q = 1}^{n_q}s_{ijq}} 
                    \end{equation}
                \item \textit{Maximise committee member suitability}. This criterion considers a product, $r_{iq}\overline{t}_{jq}x_{ijtk\ell p}$, which will be 0 unless a committee member, $i$, is assigned to a defence, $j$, that is $x_{ijtk\ell p} = 1$, and a research subject, $q$, is within the areas of expertise of a committee member, $i$, that is, $r_{iq} = 1$, and the subjects addressed in the defence, $j$, that is, $\overline{t}_{jq} = 1$, in which case the product will be 1. We want to maximise the sum of these products with this objective.
                    \begin{equation}
                        {\displaystyle \max z_3(x) = \sum_{i = 1}^{n_i} \sum_{q = 1}^{n_q} \sum_{j = 1}^{n_j}\sum_{t = 1}^{n_t}\sum_{k = 1}^{n_k}\sum_{\ell = 1}^{n_\ell}\sum_{p = 1}^{n_p} r_{iq}\overline{t}_{jq}x_{ijtk\ell p}}
                    \end{equation}
        \end{enumerate}
        \item \textit{Point of view of schedule quality}. This point of view is operationalised by the following objectives.
            \begin{enumerate}
                \item \textit{Minimise non-consecutive assignments}. The highest potential compactness value for an assignment for a committee member, $i$, is parameter $n_{v_i}$. Thus, for a committee member, $i$, the maximum potential sum of compactness values would be the product of said maximum value by the committee member's number, $j$, of assignments minus one assignment, as, logically, the committee member's first assignment cannot be scheduled within a certain time-frame after another has ended. Thus, and by weighting each assignment by the weight conferred to its committee member, the maximum sum of compactness values for the problem is represented by $\sum_{i = 1}^{n_i}\sum_{j = 1}^{n_j} u_in_{v_i}(j - 1)w_{ij}$. On the contrary, the effective sum of the compactness values, weighted by their correspondent committee member's weight, is represented by $\sum_{i = 1}^{n_i}\sum_{k = 1}^{n_k}\sum_{\ell = 1}^{n_\ell} u_{i}\overline{s}_{ik\ell}$. This criterion considers the difference between the compactness value for an ideal schedule and the effective value for the proposed schedule. We want to minimise such differences with this objective.
                    \begin{equation}\label{ob4}
                        {\displaystyle \min z_4(w, \overline{s}) = \sum_{i = 1}^{n_i}\sum_{j = 1}^{n_j} u_in_{v_i}(j - 1)w_{ij} - \sum_{i = 1}^{n_i}\sum_{k = 1}^{n_k}\sum_{\ell = 1}^{n_\ell} u_{i}\overline{s}_{ik\ell}}
                    \end{equation}
                \item \textit{Minimise the non-satisfaction of time slot preferences}. This criterion considers a product, $u_{i}(l_{ik \ell} - 1)x_{ijtk\ell p}$, which represents the penalty of assigning a committee member, $i$, to a day, $k$, at an hour, $\ell$, that is $x_{ijtk\ell p} = 1$. In this case, it takes the value of the product of the committee member's weight, $u_i$, with the penalty level assigned by the committee member to the combination of day, $k$, at an hour, $\ell$, that is, $l_{ik \ell} - 1$. The parameter, $l_{ik \ell}$, can take the value 0 if a committee member, $i$, is unavailable at a day, $k$, and an hour, $\ell$, but such an assignment would be infeasible, or a natural number, with larger values representing lower preference levels. We want to minimise the sum of such penalties with this objective.
                    \begin{equation}
                        {\displaystyle \min z_5(x) = \sum_{i = 1}^{n_i}\sum_{k = 1}^{n_k} \sum_{\ell = 1}^{n_\ell}\sum_{j = 1}^{n_j}\sum_{t = 1}^{n_t}\sum_{p = 1}^{n_p} u_{i}(l_{ik \ell} - 1)x_{ijtk\ell p}}
                    \end{equation}
                \item \textit{Minimise committee days}. This criterion takes into account a product,  $u_{i}k^{2}\overline{w}_{ik}$. The variable, $\overline{w}_{ik}$, takes the value 1 if a committee member, $i$, is assigned to committees in a number, $k$, of days, and 0 otherwise. By multiplying it by $k^2$, we keep the model linear while applying an exponential penalty representing the square of the number of days a committee member, $i$, is assigned to, promoting fairness in the committee days distribution. Moreover, the committee member's weight, $u_i$, is also considered. We want  a fair distribution of committee days between committee members and to minimise the sum of such a product with this objective. 
                    \begin{equation}
                        {\displaystyle \min z_6(\overline{w}) = \sum_{i = 1}^{n_i}\sum_{k = 1}^{n_k}u_{i}k^{2}\overline{w}_{ik}}
                    \end{equation}
                \item \textit{Minimise room changes}. The variable $\hat{s}_{ik\ell p}$ is the assigned room change penalty for a committee member, $i$, who will change rooms within a problematic time frame due to being assigned to any defence in a day, $k$, an hour, $\ell$, and a room, $p$. The penalty is incurred on the condition that the committee member, $i$, is also assigned to another defence, in a different room, $\overline{p}$, in the same day, $k$, which ends within a specific time-frame, between $\ell$ and $\ell - a_i$, with the parameter, $a_i$, representing the number of hour slots before an assignment where the committee member considers scheduling the end of another assignment as problematic for a room change. We want to minimise the sum of the product of such penalties by their incurring committee member's weight with this objective.
                    \begin{equation}\label{ob7}
                        {\displaystyle \min z_7(\hat{s}) = \sum_{i = 1}^{n_i}\sum_{j = 1}^{n_j} u_i\hat{s}_{ij}}
                    \end{equation}
            \end{enumerate}
    \end{enumerate}

\subsection{Summary}

\noindent The aim of solving scheduling problems can be analysed through the 6W's question framework. In thesis defence scheduling, we want to know Who (the committee member) is to be assigned to Which examination committee, performing What role,
to Whom (students/thesis defences), When (day and hour) and Where (the specific room).

In every MOMILP problem, the feasible region of the decision space is defined by a certain number of constraints. In our thesis defence scheduling model, we divided these constraints into three groups:

\begin{enumerate}
    \item \textit{Scheduling complete committees}. These constraints define a complete committee and ensure that every schedulable defence is assigned one.
    \item \textit{Committee composition}. These constraints ensure the eligibility of the committee members to perform their assignments.
    \item \textit{Committee member and room availability}. These constraints guarantee that committee members and rooms are available for each assignment.
\end{enumerate}

Nonetheless, not all feasible solutions are equivalent. To assess their relative standing, we identified two points of view, rendered operational by specific criteria:

\begin{enumerate}
    \item \textit{Committee assignment}. This point of view includes criteria related to the committees assigned to each defence.
    \item \textit{Schedule quality}. This point of view includes criteria for minimising the inconvenience the assignments generate for the committee members.
\end{enumerate}

Let us note that a different group of constraints is also included. Nonetheless, it aims to calculate the necessary values used to compute the various criteria, and, ergo, this group does not affect the feasible region regarding the decision variable inferred from the 6W's question framework.

\begin{figure}[h]
    \centering
    \includegraphics[scale=0.4]{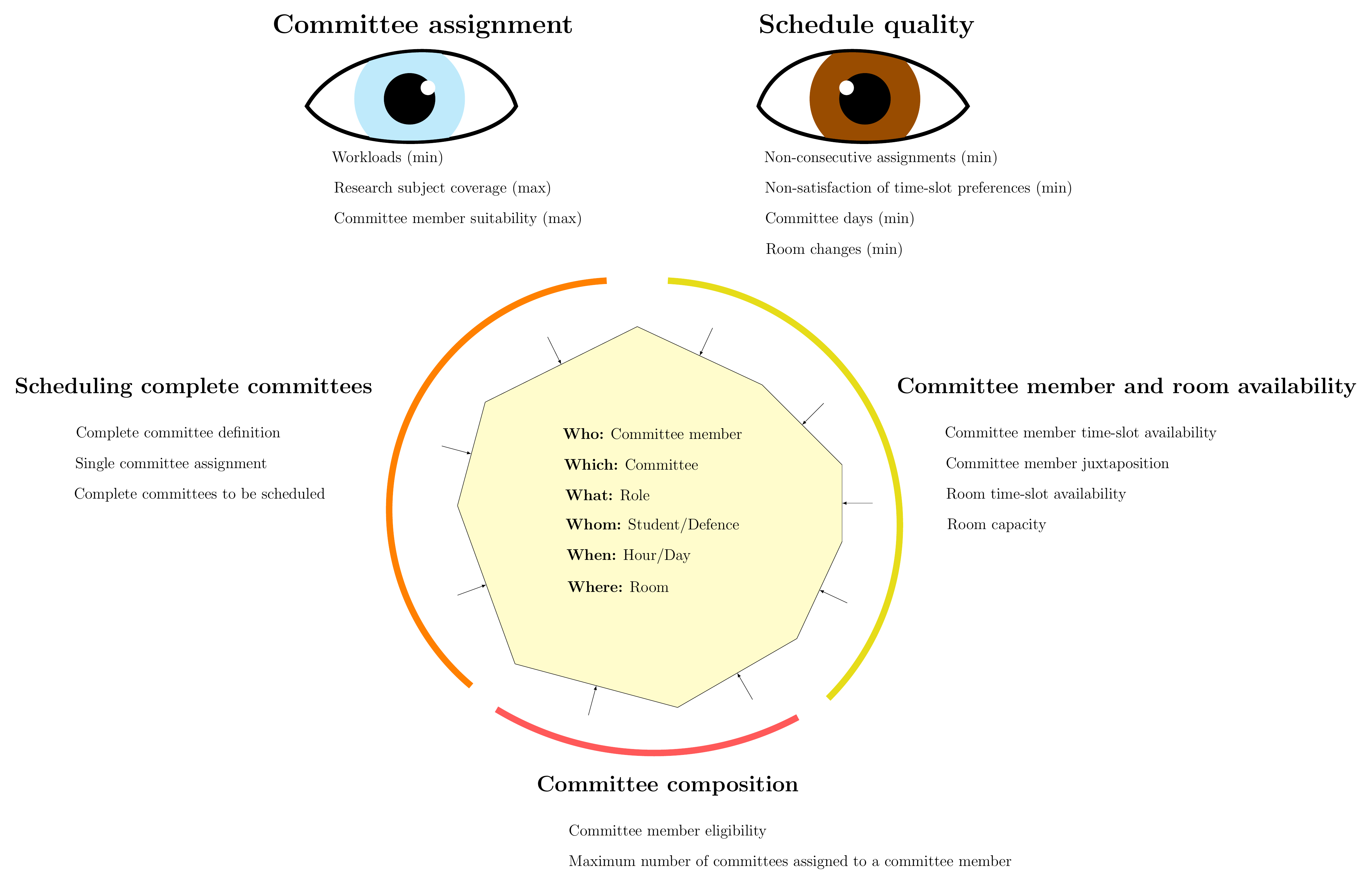}
    \caption{Thesis defence scheduling model diagram}
    \label{model sum}
\end{figure}

A representative diagram of our model is displayed in Figure \ref{model sum}.
\section{An approach to the MOMILP problem}
\noindent This section addresses the approach chosen to solve the MOMILP problem introduced in the previous section. It starts by presenting some fundamental concepts and definitions, followed by the employed algorithm and, finally, the practical aspects regarding the actual use of the algorithm.
\subsection{Some fundamental concepts, their definitions, and notation}
\noindent Consider the following MOMILP problem,

\begin{equation}
    \begin{array}{l}
        \max \;z_1(x),  \\
        \;\;\;\;\; \vdots \\
        \max \;z_i(x), \\
        \;\;\;\;\;  \vdots\\
        \max\; z_{n_z}(x),\\\\
        \text{subject to:}\\
        x \in X
    \end{array}
\end{equation}

\noindent where $x=(x_1,\ldots,x_j,\ldots,x_{n_x})$ is the vector of the decision variables, $X$ is the feasible region in the decision space, and $z_i$ is the $i-th$ linear objective function, for $i=1,\ldots,n_z$. The image of $X$ according to all the objective functions defines the feasible region, $Z$, in the criterion space.  

A fundamental concept in multi-objective optimisation is the notion of \textit{dominance}. A solution or outcome vector $z^\prime$ in the objective space dominates another solution $z^{\prime\prime}$ if and only if $z^{\prime}_i \geqslant z^{\prime\prime}_i$, for all $i=1,\ldots,n_z$, with at least one of these being a strict inequality, i.e., $z^{\prime}_i < z^{\prime\prime}_i$ for some $i$.

A feasible solution, $\bar{z} \in Z$, is said to be \textit{non-dominated} if and only if there is no other feasible solution, $z \in Z$, such that $z$ dominates $\bar{z}$. The set of all non-dominated solutions is known in the literature as the \textit{Pareto front}. The inverse image of a non-dominated solution, $\bar{x} = F^{-1}(\bar{z})$ is called an \textit{efficient solution} (in the decision space). 

Our objective is to identify a subset of the Pareto front, denoted by $N$. For this purpose, we will use a well-known scalarisation technique based on the resolution
of a sequence of $\epsilon-$constraint problems of the form. 

\begin{equation}
    \begin{array}{l}
        \max \;z_1(x),  \\
       \\
        \text{subject to:}\\
        x \in X,\\
        z_i(x) \geqslant \epsilon_i, \; i = 2, \ldots, n_z
    \end{array}
\end{equation}

\noindent where only one of the objective functions, arbitrarily chosen here as $z_1(x)$, is being maximised. Whereas the others are instead included in constraints, which set lower bounds, $\epsilon_i$, for each of the remaining objective functions, $z_i(x)$, for $i=2,\ldots,n_z$. Moreover, different non-dominated solutions are found by setting different values for the lower bounds, $\epsilon_i$.

\subsection{Algorithmic framework}

\noindent This subsection addresses the algorithmic framework. It is divided into three sequential steps: finding the number of schedulable defences, $g$, which is then set as a parameter for the subsequent steps, computing the ideal, $z^{id}$, and nadir, $z^{nad}$, points, necessary for calculating some parameters during the augmented $\epsilon-$constraint method, and, finally, the augmented $\epsilon-$constraint method itself.

In fact, the proposed method adapts the augmented $\epsilon-$constraint method, introduced in \cite{Mavrotas2009} and \cite{MavrotasFlorios2013}. This method obtains a subset of the Pareto front by iteratively increasing the lower bounds, $\epsilon_i$, for each objective, $i=2,\ldots,n_z$. However, in contrast with some $\epsilon-$constraint methods, this method guarantees that all solutions found are non-dominated.  To achieve this, instead of just using the objective function $z_1(x)$, a component related to the remaining objective functions is also included in the objective function with the help of surplus variables.

\subsubsection{Computing the ideal and the approximate nadir points}

\noindent For each $i=1,\ldots,n_z$, the problem of Equation \eqref{max objective} is solved. Each objective function is divided into two components. The first is an objective function, $z_i(x)$, for $i=1,\ldots,n_z$. The second is the sum of the remaining objective functions $z_j(x)$, for $j=1,\ldots,n_z,\;j\neq i$, multiplied by a suitable number $10^{-E}$ that ensures that, regardless of the value the second component takes, the value of the first component is the same as if we would instead maximise $z_i(x)$ separately, i.e., $\max_{x \in X} z_i(x), \; i=1,\ldots,n_z$.

\begin{equation}\label{max objective}
    \displaystyle
     z_{i}^{\rho_i\ast} = \max_{x \in X} \left\{ z_i(x) + \left(10^{-E}\right) \sum_{j = 1,\;i\neq j}^{n_z} z_{j}(x)\right\}, \;\;\; i = 1, \ldots, n_z. 
\end{equation}

Let us denote the perturbation of the objective function, $z_i(x)$, in the previous equation by,

\begin{equation}
    \rho_i = \left(10^{-E}\right) \sum_{j = 1,\;i\neq j}^{n_z} z_{j}(x), \;\; i = 1,\ldots, n_z, 
\end{equation}

\noindent where $0 \leqslant \rho_i < 1$, i.e., since $z_i(x)\in \mathbb{Z}$, then $\rho_i$ will not influence the value of $z_i(x)$. We can now define,  $z_{i}^{\ast} = z_{i}^{\rho_i\ast} - \rho_i$, for all $i = 1,\ldots, n_z$, and the ideal point, $z^{id}$, can be defined as follows:

\begin{equation}\label{ideal}
    \displaystyle
    z^{id} = \big(z_1^{\ast},\ldots,z_i^{\ast},\ldots ,z_{n_z}^{\ast}\big).
\end{equation}

Let $z_j^{i\ast}$, denote the value obtained for an objective function, $z_j(x),\;j=1,\ldots,n_z$, when maximising $z_i(x),\;i=1,\ldots,n_z$, i.e., when solving the problem of Equation \eqref{max objective}, for $z_i(x)$. Each scalar component, $z_i^{nad}$, of the approximate nadir vector, $z^{nad}$, can be defined as follows:

\begin{equation}
    \displaystyle 
    z_i^{nad} = \min_{j=1,\ldots,n_z}\left\{ z_i^{j\ast} \right\},\;\; i =1,\ldots,n_z.
\end{equation}

Accordingly, the approximate nadir vector can be stated as follows:

\begin{equation}\label{ideal}
    \displaystyle
    z^{nad} = \big(z_1^{nad},\ldots,z_i^{nad},\ldots ,z_{n_z}^{nad}\big).
\end{equation}

This process is represented in Algorithm 1, denoting the problem of Equation \eqref{max objective}, for each $z_i(x)$, as $P^{z_i}$.

\begin{algorithm}\label{Alg - maxmin}
	\caption{Compute $z^{id}$ and $z^{nad}$}
\includegraphics[scale=1.15]{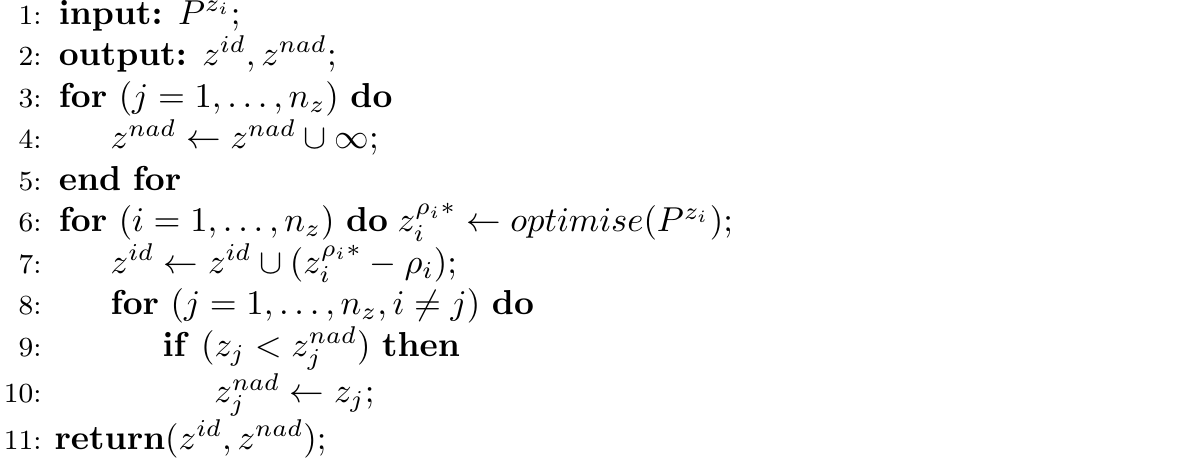}
\end{algorithm}

\subsubsection{Augmented $\epsilon-$constraint}

\noindent Finally, we can define the surplus variables, $s_i(x)$, for $x\in X$, and for each objective, $z_i(x)$, for $i=2,\ldots,n_z$, as in Equation  \eqref{equation surplus}. These are computed as the ratio of the difference between the objective, $z_i(x)$, and its corresponding value in the nadir point, $z_i^{nad}$, and the difference between the optimum value for the objective, $z_i^*$, and its corresponding value in the nadir point, $z_i^{nad}$.

\begin{equation}\label{equation surplus}
    \displaystyle s_i(x) = \frac{z_i(x) - z_i^{nad}}{z_i^{*} - z_i^{nad}}, \;\; i = 2,\ldots,n_z, \;x\in X.
\end{equation}

Now, to define the objective function used in the augmented $\epsilon-$constraint method, $z^\epsilon$, as in Equation  \eqref{obj - e const}. This objective function is divided into two components. The objective function $z_1(x)$ and a perturbation component, computed through the sum of the surplus variables, $s_i(x)$, employed to guarantee non-dominated solutions, while ensuring that the value of the objective function $z_1(x)$ is the same as if we would instead maximise $z_1(x)$ separately, i.e., $\max_{x \in X} z_1(x)$, while considering that the remaining objectives $z_j(x)$ are subject to the lower bound vector $\epsilon$, i.e., $z_j(x)\geqslant\epsilon_j$, for $j=2,\ldots,n_z$.

\begin{equation}\label{obj - e const}
    \displaystyle
    z^\epsilon =\max_{x \in X} \left\{ z_1(x) + (n_z - 0.9)^{-1} \sum_{i = 2}^{n_z} s_i(x) \right\}
\end{equation}

Let us denote the perturbation of the objective function $z_1(x)$ by

\begin{equation}
    \phi=(n_z - 0.9)^{-1} \sum_{i = 2}^{n_z} s_i(x),
\end{equation}

\noindent where, given that $0\leqslant s_i(x)\leqslant1$, it follows that $0\leqslant\phi<1$, i.e., since $z_1(x)\in \mathbb{Z}$, then $\phi$ will not influence the optimality of $z_1(x)$, considering that the remaining objective functions, $z_j(x)$, are subject to the vector of lower bounds, $\epsilon$.

To iterate between the vector of lower bounds, $\epsilon$, we first specify a parameter, $p_i$, subject to $\frac{1}{p_i}\in \mathbb{N}$, which is the percentage of the gap between the correspondent ideal, $z_i^{*}$, and nadir, $z_i^{nad}$, scalar components, which is incremented between each iteration for an objective, $z_i(x)$. Moreover, we also define a vector, $v$, of dimension $n_z - 1$, such that each scalar component, $v_i$, represents the number of times each lower bound, $\epsilon_i$, is to be incremented by its corresponding percentage, $p_i$, in the current iteration. Thus, Equation \eqref{lower bound e} defines the lower bound, $\epsilon_i$, for a given objective, $z_i(x)$, as the increment of the scalar component, $z_i^{nad}$, by a percentage, $v_ip_i$, of the difference between the scalar components, $z_i^{*}$ and $z_i^{nad}$. Note that, while all other parameters are constants, $v_i$, is updated between each iteration.  Moreover, the minimum value for the lower bound, $\epsilon_i$, must be $z_i^{nad}$, and its maximum value $z_i^{*}$, therefore $v_i \in \{0, \ldots, \frac{1}{p_i}\}$. 

\begin{equation}\label{lower bound e}
    \displaystyle
    \epsilon_i = z_i^{nad} + v_i p_i (z_i^{*} - z_i^{nad}), \; i = 2, \ldots, n_z
\end{equation}

To update $v$ between iterations, we use Algorithm 2. As input, the algorithm receives the vector $v$, the vector of percentages, $p$, and the variable $stop$, returning the updated values for $v$ and $stop$ as output. The final iteration is reached and $stop$ is set to $true$ when for all bounded objectives $z_i(x)$, $v_i = \frac{1}{p_i}$, after which the main algorithm stops and the set of the obtained non-dominated solutions, $N$, is returned. Otherwise, following an ascending order of indexes, $i = 2, \ldots, n_z$, the first $v_i < \frac{1}{p_i}$, is incremented by $1$. Furthermore, all $v_{\hat{i}},\;\hat{i} < i$, which necessarily are equal to $\frac{1}{p_i}$, are reset to $0$. This ensures that every possible vector $v$ is considered. 

\begin{algorithm}\label{Alg - iterate}
	\caption{Update $v$}
	\includegraphics[scale=1.15]{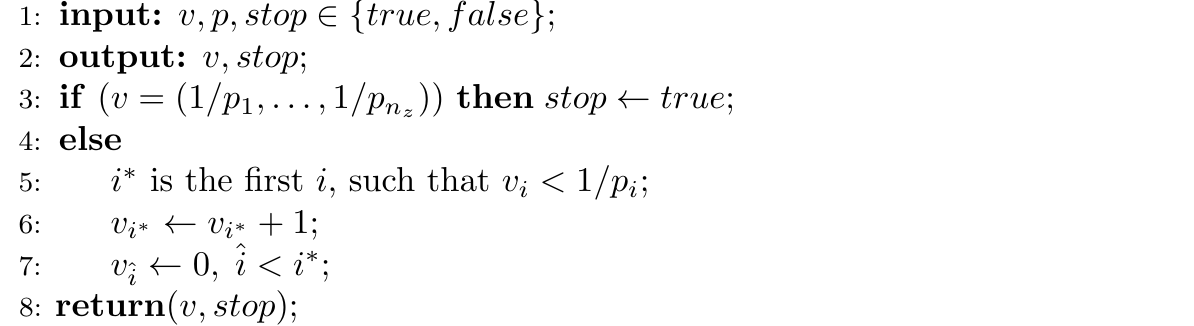}
\end{algorithm}

However, not all vectors of lower bounds, $\epsilon$, generated by vectors $v$, have the potential to obtain new solutions. Accordingly, we use Algorithm 3, which receives as input the vector of lower bounds, $\epsilon$, and the set of already obtained solutions, $N$, and assesses whether there is no already obtained solution, $z \in N$, that dominates the current lower bounds vector, $\epsilon$. If there is such a solution, $z$, then the current lower bounds vector, $\epsilon$, does not have the potential to generate a new solution, and the variable $skip$ is set to $true$. If $skip = true$, the main algorithm skips the optimisation of the problem of Equation \eqref{obj - e const} for the current $\epsilon$ (which we denote as $P^\epsilon$), improving its overall efficiency.

\begin{algorithm}
	\caption{Skip obtained solutions: $SkipSolutions()$}
\includegraphics[scale=1.15]{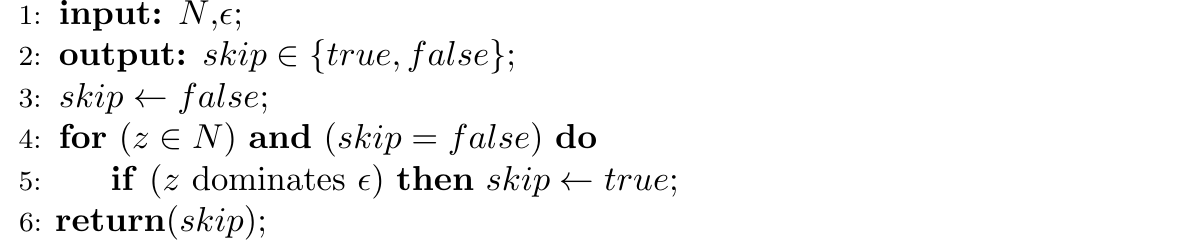}
\end{algorithm}

Conversely, some combinations can be \textit{a priori} proven to be infeasible. Similarly to the previous algorithm, Algorithm 4 tests if an iteration has the potential to enrich the pool of solutions. Accordingly, it receives as input the vector of lower bounds, $\epsilon$, and the set of lower bounds already identified as infeasible, $I$. If $\epsilon$ dominates any lower bound combination in $I$ it follows that the corresponding optimisation cannot find a feasible solution. Therefore, it can be skipped. Similarly to Algorithm 3, the variable $skip$ is set and returned as the output.

\begin{algorithm}
	\caption{Skip infeasible models: $SkipInfModels()$}
\includegraphics[scale=1.15]{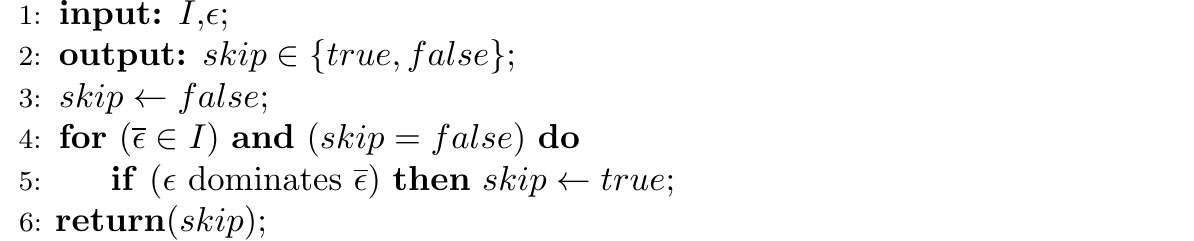}
\end{algorithm}

Finally, we can design Algorithm 5, which presents the overall augmented $\epsilon-$constraint method. As input, it receives the problem, $P^\epsilon$, the percentages vector, $p$, to be considered, and the nadir, 
$z^{nad}$, and ideal, $z_i^{*}$, points. , As output, it returns the set of obtained non-dominated solutions, $N$. 

The first step is to set up the initial maximisation, which considers $\epsilon_i = z_i^{nad}, \; i = 2, \ldots, n_z$. Additionally, the obtained solutions set, $N$, and the known infeasible lower bound combinations set, $I$, must also be defined. Necessarily, the first iteration always produces a feasible solution to be saved in the set of already obtained solutions, $N$. Nonetheless, thereupon, before every iteration, Algorithm 3 and Algorithm 4 assess if the vector of lower bounds, $\epsilon$, generated by vector $v$, has the potential to find a new solution. If it does, Problem $P^\epsilon$ is optimised, considering the current vector of lower bounds, $\epsilon$, in which case, if a new solution, $z$, is found, it is kept in the solutions set, $N$, otherwise, the lower bounds vector, $\epsilon$, is saved as infeasible in set $I$. Regardless of the present optimisation being skipped or not, Algorithm 2 determines if the last iteration has already been reached, in which case Algorithm 5 stops and returns the obtained solutions set, $N$. Otherwise, Algorithm 2 updates $v$ for the next iteration.

\begin{algorithm}[h]
	\caption{Augmented $\epsilon-$ constraint}
	\includegraphics[scale=1.15]{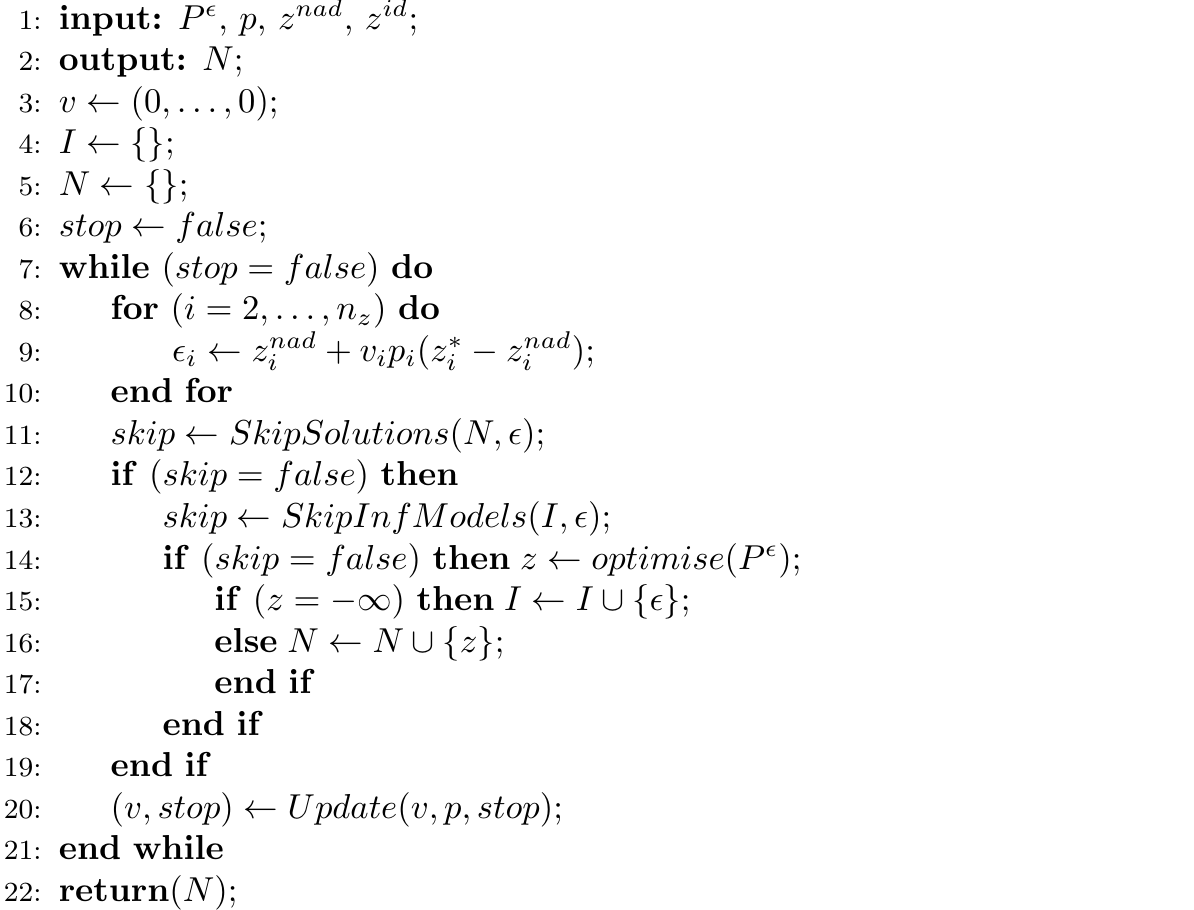}
\end{algorithm} 

\subsubsection{Finding the number of schedulable defences}
\noindent Before computing $z^{id}$, $z^{nad}$, and initialising Algorithm 5, we need to set the value for the number of defences that can be scheduled in any given instance, $g$. For such a purpose, we solve the alternative problem of Equation \eqref{ob8}, which is similar to Problem $P^{\epsilon}$, but without considering Constraint \eqref{const - defences scheduled once - eq 5}, which sets the number of defences that are to be scheduled, Constraints \eqref{const - res perc - eq 1}-\eqref{const - room - eq 4}, which define the values for the objectives, and lastly, the objectives themselves \eqref{ob1}-\eqref{ob7}.

Conversely, we instead include the objective function of Equation \eqref{ob8}, which maximises the number of scheduled complete committees (thesis defences), computed as the sum of a variable, $y_{jk\ell p}$, which takes the value 1 if a defence, $j$, is scheduled at a day, $k$, an hour, $\ell$, and a room, $p$. 

\begin{equation}\label{ob8}
    \displaystyle g = \max z_g(y) = \sum_{j = 1}^{n_j}\sum_{k = 1}^{n_k}\sum_{\ell = 1}^{n_\ell}\sum_{p = 1}^{n_p} y_{jk\ell p}
\end{equation}

\begin{figure}[htbp!]
    \centering
   \begin{tikzpicture}[node distance = 1cm, auto]

    \node[block4] (start) {Start};
    \node[bigbox1, below = of start, xshift=1cm, yshift=0.35cm] (1st stage) {First Stage};
    \node[bigbox, below = of 1st stage, yshift=0.7cm] (1st stage) {Second Stage};
    \node[block5, below = of start] (G) {$g \leftarrow z_g(y)$};
    \node[block5, below = of G] (Alg1) {Algorithm 1};
    \node[block5, below = of Alg1] (Alg5) {Algorithm 5};
    \node[block4, below = of Alg5] (end) {End};
    
    \path [line] (start) -- (G); 
    \path [line] (G) -- (Alg1); 
    \path [line] (Alg1) -- (Alg5); 
    \path [line] (Alg5) -- (end); 
    \end{tikzpicture}
\caption{Full procedure diagram}\label{fullprocedure}
\end{figure}
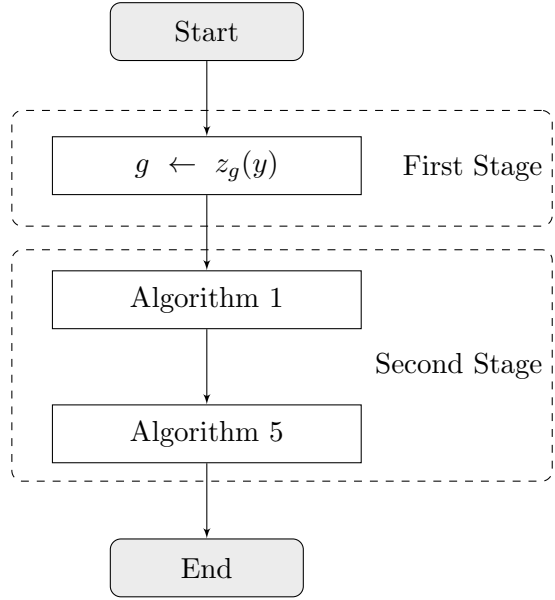

To sum up, the first stage of the procedure is to find the maximum number of thesis defences that can be scheduled for a given instance and set that value as a parameter in the following steps. For the second stage, we compute the ideal point, $z^{id}$, and the approximate nadir point, $z^{nad}$ through Algorithm 1. Finally, we have all the necessary parameters to initialise the augmented $\epsilon-$constraint method, Algorithm 5.  A schematic representation of the whole two-stage procedure is presented in Figure \ref{fullprocedure}.

\section{Designing an instance generator for thesis defence scheduling}

\noindent This section addresses the design of the instance generator for thesis defence scheduling problems that we propose. It starts by presenting the different types of instances that were considered and references the different additional necessary inputs. Moreover, the specific procedures for the random instance generation are also presented, considering the parameters generated through simple random choices and the availability parameters, defined based on conditional probabilities.

\subsection{Types of instances and other inputs}\label{sec inputs}

\noindent To test the proposed method (see Section \ref{sec comp exp}), a set of 96 instances, denoted by $p({n_i.n_j.n_t.n_k.n_\ell. n_p.n_q})$, was generated. Let us point out that the following parameters are identical for all of these instances: 

\begin{itemize}[label={--}]
    \item $n_t = 3$, which is the defined number of roles;
    \item $n_k = 15$, which is the defined number of days;
    \item $n_\ell = 16$, which is the defined number of hour slots in a day;
    \item $n_q = 15$, which is the defined number of research subjects.
\end{itemize}

Moreover, the instances were divided into six different types by varying the number of committee members, $n_i$, the number of defences, $n_j$, and the number of rooms, $n_p$:

\begin{enumerate}
    \item Instances of type, $p({25.20.3.15.16.3.15})$: These instances consider 25 committee members ($n_i = 25$), 20 defences ($n_j = 20$) and 3 rooms ($n_p = 3$), instances (1)-(16).
    \item  Instances of type, $p({25.20.3.15.16.4.15})$: These instances consider 25 committee members ($n_i = 25$), 20 defences ($n_j = 20$) and 4 rooms ($n_p = 4$), instances (17)-(32).
    \item  Instances of type, $p({38.30.3.15.16.3.15})$: These instances consider 38 committee members ($n_i = 38$), 30 defences ($n_j = 30$) and 3 rooms ($n_p = 3$), instances (33)-(48).
    \item  Instances of type, $p({38.30.3.15.16.4.15})$: These instances consider 38 committee members ($n_i = 38$), 30 defences ($n_j = 30$) and 4 rooms ($n_p = 4$), instances (49)-(64).
    \item  Instances of type, $p({50.40.3.15.16.3.15})$: These instances consider 50 committee members ($n_i = 50$), 40 defences ($n_j = 40$) and 3 rooms ($n_p = 3$), instances (65)-(80).
    \item  Instances of type, $p({50.40.3.15.16.4.15})$: These instances consider 50 committee members ($n_i = 50$), 40 defences ($n_j = 40$) and 4 rooms ($n_p = 4$), instances (81)-(96).
\end{enumerate}

As for the remaining parameters, some were considered the same for every instance, specifically:

\begin{itemize}[label={--}]
    \item $d=2$, which is the duration of a thesis defence, in the number of time slots;
    \item $p(u_i = 1) = 0.7, \;p(u_i = 2) = 0.3$, which are the probabilities of a committee member $i$ to be assigned to a certain individual weight $u_i$;
    \item $c_i = 0.5n_i$, which is the maximum number of allowed defences \textit{per} committee member $i$;
    \item $\sum_{q = 1}^{n_q} r_{iq} = 3$, which is the number of research subjects $q$ for committee members $i$;
    \item $\sum_{q = 1}^{n_q} t_{jq} = 3$, which is the number of research subjects $q$ for defences $j$.
\end{itemize}

Additionally, for each type of instance, we generated sixteen different randomised instances while varying some data inputs corresponding to the remaining parameters of our model: 

\begin{enumerate}
    \item $t=1$ and $t=2$, or $t=2$: These are the fixed roles in $e_{ijt}$. If a role, $t$, is fixed, it means that, for each defence, $j$, there is only one eligible committee member, $i$, it can be assigned to. Regardless of a role being fixed, there is an overall set of committee members that can be assigned to that role. For the instances with 25 committee members ($n_i = 25$), we defined that 9 committee members are to be randomly selected as eligible for role $t = 1$ and 13 for role $t = 2$. For the instances with 38 committee members ($n_i = 38$), 12 committee members are selected as eligible for role $t = 1$ and 19 for role $t = 2$. Lastly, for the instances with 50 committee members ($n_i = 50$), 15 committee members are selected as eligible for role $t = 1$ and 25 for role $t = 2$. Let us note that it is possible for a committee member to be part of the randomly selected sample for different roles.  Finally, all committee members are eligible for role $t=3$.
    \item $p(l_{ik\ell}=0)=0.78, \; p(l_{ik\ell}=0)=0.82$, or $p(l_{ik\ell}=0)=0.86$: These are the unavailability percentages, $p(l_{ik\ell}=0)$, for the committee member time slot preference and availability parameter, $l_{ik\ell}$, which are defined as the percentage of $l_{ik\ell}=0$. The percentages for the instances with 2 fixed roles were 0.78 or 0.82. Conversely, for the instances with a single fixed role, they were 0.82 or 0.86. This differentiation between instances with 1 or 2 fixed roles was necessary due to increased computational complexity.
    \item $p(m_{k\ell p}=0)=0.8$ or $p(m_{k\ell p}=0)=0.86$, for $m_{k\ell p}$: These are the unavailability percentages, $p(m_{k\ell p}=0)$, for the room availability parameter, $m_{k\ell p}$, which are defined as the percentage of $m_{k\ell p}=0$.
    \item $p(v_{i\ell}=[1])=0.7,\;p(v_{i\ell}=[2,1])=0.3$ or $p(v_{i\ell}=[1])=0.8,\;p(v_{i\ell}=[2,1])=0.2$: These are the probabilities of a committee member, $i$, being assigned values [1] or [2,1] for the compactness preference parameter, $v_{i\ell}$, which will affect the big-$M$ upper bounds, $n_{v_i}$.
    \item $p(h_{i\ell}=[1])=0.7,\;p(h_{i\ell}=[2,1])=0.3$ or $p(h_{i\ell}=[1])=0.8,\;p(h_{i\ell}=[2,1])=0.2$: These are the probabilities of a committee member, $i$, being assigned values [1] or [2,1] for the room change penalty parameter, $h_{i\ell}$, which will affect the big-$M$ upper bounds, $n_{h_i}$.
\end{enumerate}

\begin{figure}[h]
    \centering
    \includegraphics[width=8cm]{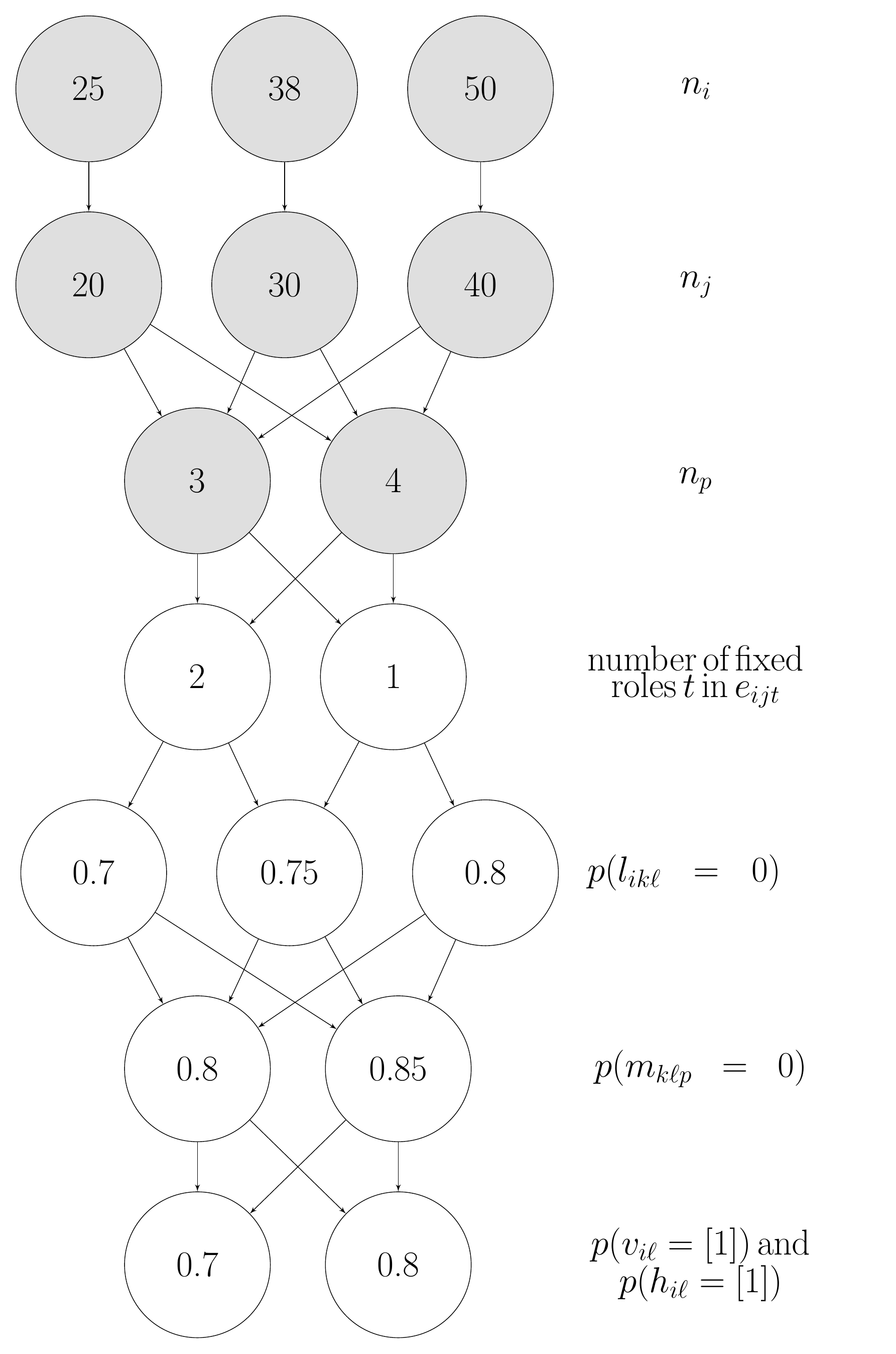}
    \caption{Generated instances diagram}
    \label{gene insts diagram}
\end{figure}

A diagram summarising this instance generation procedure is presented in Figure \ref{gene insts diagram}.

\subsection{Random choice parameters}
\noindent Regarding the different parameters that are randomly generated, some of them are obtained through simple random choices, which consider the probabilities presented in Subsection \ref{sec inputs}:
\begin{enumerate}
    \item $u_i$: For every committee member, $i$, the value 1 or 2 is selected.
    \item $v_{i\ell}$: For every committee member, $i$, the value $[1]$ or $[2, 1]$ is selected.
    \item $h_{i\ell}$: For every committee member, $i$, the value $[1]$ or $[2, 1]$ is selected.
    \item $r_{iq}$: For every committee member, $i$, a sample of 3 research subjects, $q$, out of the 15 is randomly selected, with a uniform probability distribution, and the corresponding parameter value is set to 1.
    \item $t_{jq}$: For every defence, $j$, a sample of 3 research subjects, $q$, out of the 15 is randomly selected, with a uniform probability distribution, and the corresponding parameter value is set to 1.
    \item $e_{ijt}$: For fixed roles, $t$, and for each defence, $j$, out of the randomly selected eligible committee members sample, a committee member, $i$, is selected with a uniform probability, and the corresponding value for $e_{ijt}$ is set to 1, and all other values of $e_{ijt}$ for that defence, $j$, and role, $t$, are set to 0. We denote that when two or more roles, $t$, are fixed, the same committee member, $i$, cannot be assigned two roles for a given defence, $j$. For the instances where the role $t = 1$ is not fixed, for all committee members, $i$, that have been randomly selected as eligible and all defences, $j$, the parameter, $e_{ij1}$, is set to 1. Conversely, for role $t = 3$, the parameter, $e_{ij3}$, is always set to 1.
\end{enumerate}

\subsection{Availability parameters generation for committee members and rooms}

\noindent In real-world thesis defence scheduling applications, the availability periods for committee members and rooms usually occur in blocks, between lectures, before the first lecture of a given day, or after all the daily assignments. Algorithm 6 was designed to replicate such behaviour in the generation of the availability of committee members, $l_{ik\ell}$, and rooms, $m_{k\ell p}$. 

The algorithm can be described as a Markov chain, which will help estimate the probability distributions of the different possible values for the availability parameters, $l_{ik\ell}$ and $m_{k\ell p}$.

This subsection starts by addressing the inputs for the availability generation algorithm. Then it describes the algorithm itself. Subsequently, we present some fundamental concepts regarding Markov chains and explain why this algorithm can be classified as one. Moreover, this definition is used to estimate the probability distribution for the different values the availability parameters can take.

For simplification, when we refer to an individual $\tau$, it can represent a committee member or a room.

\subsubsection{Inputs for the generation of the availability parameters}
\noindent To generate the availability parameters in a manner which represents their real-world behaviour, we defined the probability of an individual, $\tau$, to have a certain availability status, $\alpha$, at any day, $k$, and hour, $\ell$, (time slot ($k,\ell$)), as conditional on its status, $\overline{\alpha}$, in a previous time slot, ($k, \overline{\ell}$), for $\overline{\ell} = \ell - 1$.
Let us note that, when we mention a conditional probability, $p(\alpha|\overline{\alpha})$, what we are referring to is the probability of having $l_{ik\ell} = \alpha$ or $m_{k\ell p}=\alpha$, given that $l_{ik\overline{\ell}} = \overline{\alpha}$ or $m_{k\overline{\ell} p}=\overline{\alpha}$.

Algorithm 6 receives the following inputs:

\begin{enumerate}
    \item $n_i$ or $n_p$: This is the number of committee members, $i$, or rooms, $p$.
    \item $n_k$: This is the number of days, $k$.
    \item $n_\ell$: This is the number of hour slots, $\ell$, in a day.
    \item $d$: This is the duration of the defences.
    \item $\Delta$: This is the duration of the initial warm-up period for each day, $k$, within which the generated parameters will be disregarded.
    \item $p(\alpha|\alpha)$: These are the probabilities of an $l_{ik\ell} = \alpha$ or $m_{k\ell p} = \alpha$ to remain unchanged between $\ell$ and $\overline{\ell}, \; \ell = \overline{\ell}  + 1$.
    \item $p(\alpha|\overline{\alpha})$: These are the probabilities of an $l_{ik\ell} = \alpha$ or $m_{k\ell p} = \alpha$ changing from a state, $\overline{\alpha}$, to another, $\alpha$, between $\ell$ and $\overline{\ell}, \; \ell = \overline{\ell}  + 1$. This is not an input \textit{per se}, but computed through Equation \eqref{equation conditional change}, based on the input values, $p(\alpha|\alpha)$.
\end{enumerate}

Computing the remaining probabilities, $p(\alpha|\overline{\alpha})$, through Equation \eqref{equation conditional change} promotes the proportionality between $p(\alpha|\overline{\alpha})$ and all other $p(\hat{\alpha}|\overline{\alpha})$, based on their conditional probabilities, $p(\alpha|\alpha)$. Moreover, this equation also guarantees that the sum of these probabilities is always equal to 1.

\begin{equation}\label{equation conditional change}
\displaystyle
    p(\alpha|\overline{\alpha}) = p(\alpha|\alpha)\left(\sum_{\hat{\alpha}=0}^{n_\alpha} p(\hat{\alpha}|\hat{\alpha})\right)^{-1} (1 - p(\overline{\alpha}|\overline{\alpha})), \; \alpha, \hat{\alpha} \neq \overline{\alpha} 
\end{equation}

For the committee member availability parameter, $l_{ik\ell}$, three values were considered for assignment, $l_{ik\ell} = 0$, representing unavailability, $l_{ik\ell} = 1$, representing preferred time slots, and, $l_{ik\ell} = 2$, representing less preferred time slots. The room availability parameter, $m_{k\ell p}$, was defined as binary, hence, the possible values are, $m_{k\ell p}=0$, representing unavailability, and, $m_{k\ell p}=1$, representing availability. The inputs, $p(\alpha|\alpha)$, per each defined unavailability percentage, are presented in Table \ref{gen l} for the committee member availability parameter, $p(l_{ik\ell}=0)$, and in Table \ref{gen m} for the room availability parameter, $p(m_{k\ell p}=0)$.

\begin{table}[h]
   \caption{Conditional probabilities $p(\alpha|\alpha)$ \textit{per} generated unavailability percentage $p(l_{ik\ell}=0)$ for parameter $l_{ik\ell}$}
    \centering
    \begin{tabular}{|c|ccc|}
    \hline
         &\multicolumn{3}{c|}{$p(l_{ik\ell}=0)$}\\
        \hline
        $p(\alpha|\alpha)$& 0.78& 0.82&0.86 \\
        \hline
        $p(0|0)$ & 0.95  & 0.95 & 0.95\\
        $p(1|1)$ & 0.7 & 0.63 & 0.55\\
        $p(2|2)$ & 0.7 & 0.63 & 0.55\\
        \hline
    \end{tabular}
    \label{gen l}
\end{table}

\begin{table}[h]
   \caption{Conditional probabilities $p(\alpha|\alpha)$ \textit{per} generated unavailability percentage $p(m_{k\ell p}=0)$ for parameter $m_{k\ell p}$}
    \centering
    \begin{tabular}{|c|cc|}
    \hline
         &\multicolumn{2}{c|}{$p(m_{k\ell p}=0)$}\\
        \hline
        $p(\alpha|\alpha)$& 0.8& 0.86 \\
        \hline
        $p(0|0)$ & 0.95  & 0.95 \\
        $p(1|1)$ & 0.7 & 0.8 \\
        \hline
    \end{tabular}
    \label{gen m}
\end{table}

\subsubsection{Availability generation for committee members and rooms algorithm}
\noindent The first step of the algorithm is to define how the values for the first hour slot, $\ell=1$, for an individual, $\tau$, and day, $k$, $l_{ik1}$ and $m_{k1 p}$, are generated. We determined that, for each parameter regarding $\ell = 1$, the conditional probabilities of it being assigned a status $\alpha$, $p(\alpha|\overline{\alpha})$, are to consider $l_{ik0} = 0$ or $m_{k0 p} = 0$, i.e., $\overline{\alpha}=0$. Nonetheless, this creates a biased availability distribution for the first hour slots for each day. However, we want to ensure that the probability, $p(\alpha)$, of an hour slot, $\ell$, being assigned availability status, $\alpha$, is independent of the hour slot, $\ell$. Thus, for every committee member or room, and day, we generate an excess of parameter values  to uniformise their distribution. Then we disregard the initial excess values, i.e., we consider a warm-up period, $\Delta$. In our experiments, we used a $\Delta = 40$ \textit{per} committee member or room, and day, which proved sufficient to eliminate the initial value's effect.

Then, we generate each availability parameter value in sequential order, based on the presented conditional probabilities. 

Nonetheless, there is an exception where these input probabilities do not apply. It occurs after an individual, $\tau$, which was available to start an assignment in a time slot, ($k, \overline{\ell}$), stops being available to do so at the next time slot, ($k, \ell$), for $\ell = \overline{\ell}+1$. Let us note that for an individual $\tau$ to be able to start an assignment at time slot ($k, \overline{\ell}$), its effective available time must extend up until ($k, \hat{\ell}$), for $\hat{\ell} = \ell + d$, as, otherwise, the individual, $\tau$, could not be present for the whole duration, $d$, of a defence. Thus, we defined that when a change from availability to unavailability occurs, the individual, $\tau$, must be unavailable for at least a long enough time, such that if they were scheduled for their last available time slot before the unavailability block, and for the first available one after its occurrence, the individual, $\tau$, would have at least one time slot of effectively unassigned time. Therefore, whenever such a change occurs, $d-1$ unavailable time slots are automatically added, ensuring the rule above is respected. 

\begin{algorithm}[h]
	\caption{Generate availability parameters}
\includegraphics[scale=1.15]{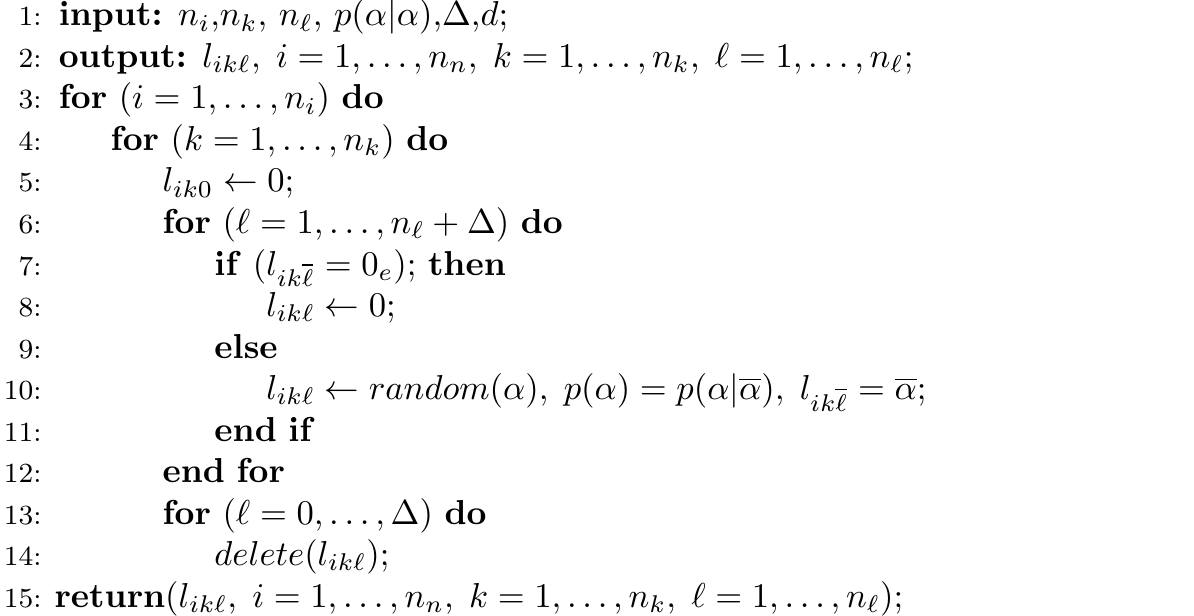}
\end{algorithm} 

An illustrative example of the time slot preference and availability generated by this method is presented in Figure \ref{Availability diagram}. Let us note that, since we considered $d = 2$, the exceptional addition was of a single time slot (i.e., $d - 1$). For example, if we had used $d=4$ instead, three time slots would have been added.

\begin{figure}[h]
\centering
\begin{tikzpicture}
 \node[blockns] (1) {non-success\\$l_{ik1}\leftarrow2$};
 \node[blockns, right = of 1, xshift=-1cm] (2) {non-success\\$l_{ik2}\leftarrow2$};
 \node[blocks, right = of 2, xshift=-1cm] (3) {success\\$l_{ik3}\leftarrow0$};
 \node[block0, right = of 3, xshift=-1cm] (4) {exceptional addition\\$l_{ik4}\leftarrow0$};
 \node[blockns, right = of 4, xshift=-1cm] (5) {non-success\\$l_{ik5}\leftarrow0$};
 \node[blocks, right = of 5, xshift=-1cm] (6) {success\\$l_{ik6}\leftarrow1$};
 \node[blockns, right = of 6, xshift=-1cm] (7) {non-success\\$l_{ik7}\leftarrow1$};
\end{tikzpicture}
   
    \caption{Committee member time slot preference and availability parameter $l_{ik\ell}$ generation illustrative diagram}
    \label{Availability diagram}
\end{figure}

\subsubsection{Some fundamental concepts regarding Markov chains}

\noindent Our availability generation algorithm can be defined as a Markov chain. Nonetheless,  some fundamental concepts must be clarified before this representation can be addressed.

A Markov chain is a type of stochastic process, with the distinguishing characteristic that each of its states, $\alpha$, is part of a set of discrete events and that the probability of each state, $\alpha$, to occur at time $\ell$, depends only on the previous state, $\overline{\alpha}$, occurring at time $\ell-1$. 

A square transition matrix $T$ $i\times j$ can also represent such a system. Each entry of transition matrix $T$ represents the probability $p(\alpha_j|{\alpha}_i)$ of the next state being $\alpha_j$ given that the previous state was ${\alpha}_i$. Moreover, each entry must be within 0 and 1, as they represent probabilities, and the sum of each row must be 1, to correctly represent the total probability of a given set. An example of such a matrix $T$ is displayed in Figure \ref{Transition M}.

\begin{figure}[h!]
    \centering
    \begin{tabular}{|c|c|c|c|c|}
        \hline
        $p(\alpha_1|{\alpha}_1)$ & \ldots &$p(\alpha_j|{\alpha}_1)$&\ldots&$p(\alpha_{n_j}|{\alpha}_{1})$\\
        \hline
        \multicolumn{5}{|c|}{\ldots}\\
         \hline
         $p(\alpha_1|{\alpha}_i)$ & \ldots &$p(\alpha_j|{\alpha}_i)$&\ldots&$p(\alpha_{n_j}|{\alpha}_{i})$\\
        \hline
        \multicolumn{5}{|c|}{\ldots}\\
        \hline
         $p({\alpha}_1|\alpha_{n_i})$ & \ldots &$p({\alpha}_j|\alpha_{n_i})$&\ldots&$p(\alpha_{n_j}|{\alpha}_{n_i})$\\
         \hline
    \end{tabular}
    \caption{An example of a transition matrix $T$}
    \label{Transition M}
\end{figure}

Accordingly, the power $T^\ell$ computes the probability of each state $\alpha_j$ to occur after $\ell$ repetitions, given that the initial state was $\alpha_i$. Moreover, a transition matrix $T$ is said to be regular if, after a certain number of repetitions, each of its columns stabilises at a certain value. Thus, if a transition matrix $T$ is regular, there is a vector $V$, such that, after a sufficiently large number of experiments, $\ell$, and for any probability vector $\hat{p}$, the following condition is verified,

\begin{equation*}
    \displaystyle \hat{p}\cdot T^{\ell} \approx V.
\end{equation*}

This suggests that after a certain number of experiments, regardless of the initial conditions, a regular Markov chain converges to a steady-state, with each possible state, $\alpha$, occurring at a certain probability $V_\alpha$. Furthermore, a transition matrix $T$ is known to be regular, if, after any number of repetitions, $\ell$, there is a $T^{\ell}$ such that $p(\alpha_j|\alpha_i)>0$, for all $ j = 1,\ldots, n_j, \;i=1\ldots,n_i$.

\subsubsection{Representation of Algorithm 6 as a Markov chain and probability distribution estimation}
\noindent Our availability generation algorithm can be described as a Markov chain with $n_\alpha + d$ possible states, where $n_\alpha$ is the number of different availability states and $d$ is the duration of the defences. 

For an availability state, $\alpha=1,\ldots,n_\alpha$, occurring at time $\ell-1$, the following state, i.e., the one occurring at time $\ell$, can remain unchanged, with a probability $p(\alpha|\alpha)$, it can change to another availability state, $\overline{\alpha}=1,\ldots,n_\alpha,\;\overline{\alpha}\neq \alpha$, with probability $p(\overline{\alpha}|\alpha)$ or it can start an unavailability block, with probability $p(0|\alpha)$. 

When an unavailability block starts, we know that at least $d$ zeros are to be added. Nonetheless, in a Markov chain, the probabilities must only depend on the previous state. Thus, before we have the unavailability state, $\alpha=0$, there must be $d-1$ different states, $0_{e_{\hat{i}}},\;\hat{i}=1,\ldots,d-1$, that lead to the generation of a 0, with probability 1, i.e., for $\hat{i}=1,\ldots,d-2$, the only probability is $p(0_{e_{\hat{i+1}}}|0_{e_{\hat{i}}})=1$ and for $\hat{i}=d-1$ it is $p(0|0_{e_{{d-1}}})=1$.

Finally, when the last exceptional 0 is added, the state $\alpha=0$ functions similarly to the availability states. Specifically, it can repeat itself, with probability $p(0|0)$, or it can generate an availability state, $\alpha=1,\ldots,n_\alpha$, with probability $p(\alpha|0)$. Let us note that, unlike the availability states, the unavailability state, $\alpha=0$, cannot be followed by a state $0_{e_1}$.

\begin{figure}[h!]
    \centering
    \includegraphics[scale=0.35]{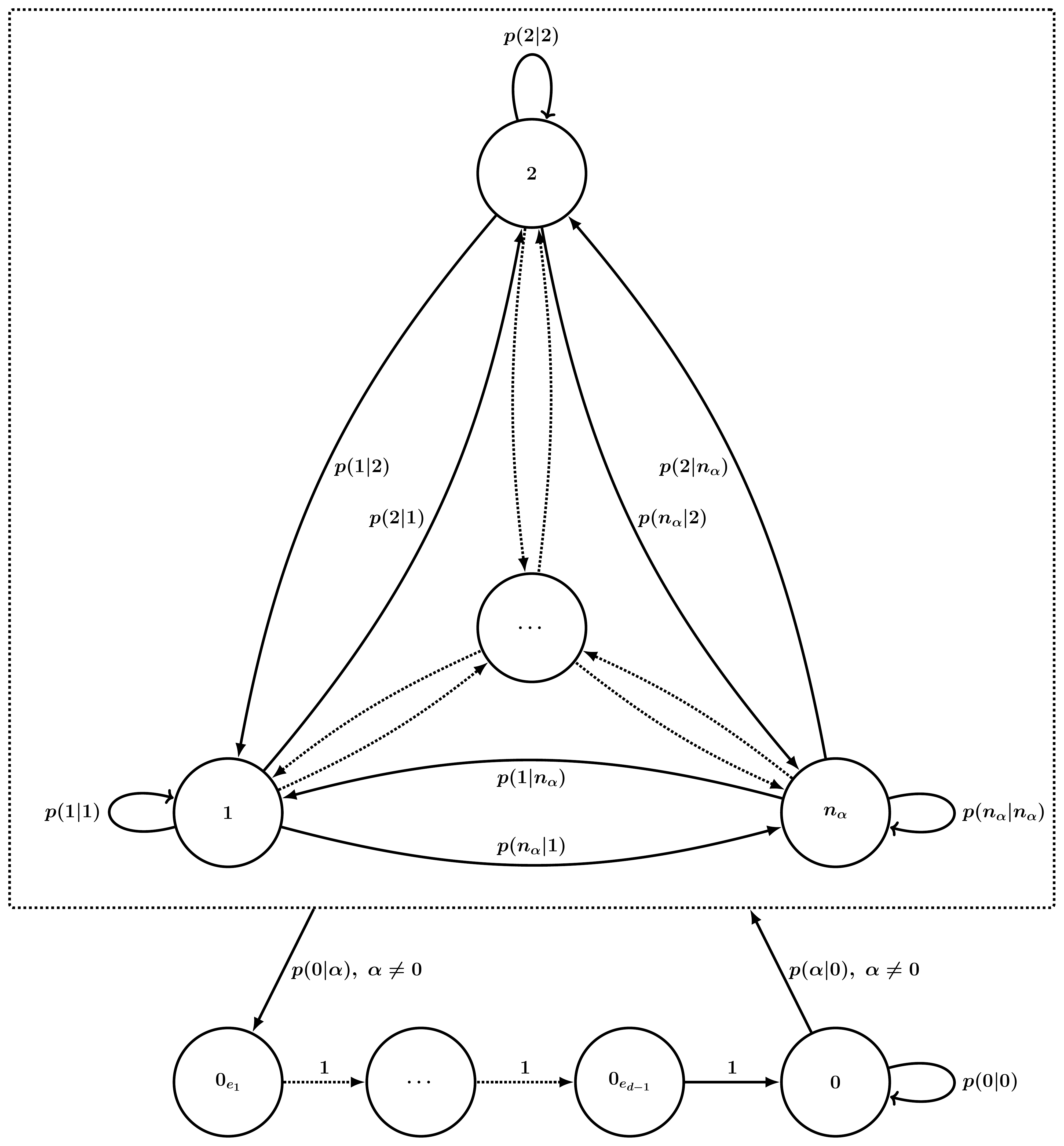}
    \caption{Markov chain representation of the availability generation algorithm}
    \label{fig:mchain}
\end{figure}

The Markov chain that represents the availability generation algorithm is displayed in Figure \ref{fig:mchain}.

Accordingly, we can use the properties of each corresponding transition matrix, $T$, to estimate the probability distribution of the parameters generated through this method, $l_{ik\ell}$ and $m_{k\ell p}$. 

Considering the warm-up period, $\Delta=40$, and, as an example, the generation of $l_{ik\ell}$ with an unavailability percentage $p(l_{ik\ell}=0)=0.78$:

The transition matrix, $T$, is:

\begin{table*}[h!]
    \centering
    \begin{tabular}{|c|c|c|c|c|}
    \hline
        	$p(\alpha_j|\alpha_i)$& $0_{e_1}$ &	0 &	1 &	2\\
        	\hline
$0_{e_1}$ &	0	&1&	0&	0\\
\hline
0&	0&	0.95&	0.025&	0.025\\
\hline
1&	0.1728&	0&	0.7	&0.1272\\
\hline
2&	0.1728	&0&	0.1272&	0.7\\
\hline
    \end{tabular}
\end{table*}

The transition matrix, following 39 repetitions, $T^{39}$, is:

\begin{table*}[h!]
    \centering
    \begin{tabular}{|c|c|c|c|c|}
    \hline
        	$p(\alpha_j|\alpha_i)$& $0_{e_1}$ &	0 &	1 &	2\\
        	\hline
$0_{e_1}$ &	0.0373	&0.7466&0.1080&	0.1080\\
\hline
0&	0.0373&0.7466&0.1080&0.1080\\
\hline
1&	0.0373&	0.7466&0.1080&0.1080\\
\hline
2&	0.0373	&0.7466&	0.1080&0.1080\\
\hline
    \end{tabular}
\end{table*}

And the transition matrix, following 40 repetitions, $T^{40}$, is:

\begin{table*}[h!]
    \centering
    \begin{tabular}{|c|c|c|c|c|}
    \hline
        	$p(\alpha_j|\alpha_i)$& $0_{e_1}$ &	0 &	1 &	2\\
        	\hline
$0_{e_1}$ &	0.0373	&0.7466&0.1080&	0.1080\\
\hline
0&	0.0373&0.7466&0.1080&0.1080\\
\hline
1&	0.0373&	0.7466&0.1080&0.1080\\
\hline
2&	0.0373	&0.7466&	0.1080&0.1080\\
\hline
    \end{tabular}
\end{table*}

Consequently, $p(l_{ik\ell}=0)=0.0373+0.7466\approx0.78$ and $p(l_{ik\ell}=1)=p(l_{ik\ell}=2)\approx0.11$. Let us note that, at a $\ell=\Delta=40$, the matrixes were not yet fully stationary, with differences between $T^{39}$ and $T^{40}$ of order $10^{-5}$. Nonetheless, we considered these to be small enough to conduct our computational experiments.

\section{Computational experiments, results, and some comments}\label{sec comp exp}

\noindent This section addresses the computational experiments made on the generated instances. It starts by specifying some essential practical aspects, and then the analysis of the results of the computational experiments is presented.

\subsection{Practical aspects}

\noindent There are several practical aspects we need to consider regarding hardware and software, iteration time limits, parameters of the augmented $\epsilon-$constraint method, and the display of the computational experiments.

\subsubsection{Hardware and software}
\noindent Regarding both hardware and software characteristics:

\begin{enumerate}
    \item CPU: 11th Gen Intel(R) Core(TM) i5-1135G7 @ 2.40GHz   2.42 GHz.
    \item RAM: 15.8GB.
    \item Implementation of the algorithms: Python 3.10.
    \item Solver: Gurobi 9.5.0.
\end{enumerate}

\subsubsection{Time limits}
\noindent Regarding the time limits set for each step of the procedure:

\begin{enumerate}
    \item Finding parameter $g$: 30 minutes.
    \item Algorithm 1: 2 hours, equally divided between the seven objectives.
    \item Algorithm 5: For each iteration, 12 hours minus the time used for the previous steps and iterations in the procedure, divided by the remaining number of iterations.
\end{enumerate}

\subsubsection{Parameters of the algorithms}
\noindent Regarding the parameters of the augmented $\epsilon-$constraint method and other algorithms:

\begin{enumerate}
    \item Objective to be fully considered in $z^\epsilon$: $z_1$.
    \item Bounded objectives: $z_3$, $z_4$.
    \item Number of equally spaced bounds between, and including, $z^{nad}_i$ and $z^{*}_i$: 10.
    \item Continuous objective $z_2$: While $z_2$ is defined as continuous, it was rounded up to the nearest integer, as, otherwise, Algorithm 1 could not be used to assess its optimum accurately.
\end{enumerate}

\subsubsection{Computational experiments display}
\noindent Regarding the computational experiments, displayed in the Appendix, in Tables \ref{small}, \ref{medium} and \ref{large}, for small ($n_j = 20$), medium ($n_j = 30$) and large ($n_j = 40$) instances, respectively:

\begin{enumerate}
    \item Table row: Presents the number of the instance, the identification of the type of instance, the corresponding input data and the outputs.
    \item Types of output: Presents the number of non-dominated solutions found, $|N|$, the number of infeasible iterations, $|I|$, the number of skipped non-dominated solutions, $skip^N$, the number of skipped infeasible solutions, $skip^I$, the number of solutions found that were not proven as non-dominated due to a time limit being reached, $time^N$, the number of iterations that were stopped before a solution was found due to a time limit being reached, $time^I$, the number of defences that can be scheduled, $g$, and the CPU time required.
\end{enumerate}

\subsection{Computational experiments}

\noindent This subsection addresses the analysis of the computational experiments. For such a purpose, several aspects must be  taken into account. Specifically, the computational performance, number of schedulable defences, and the type of iteration distribution. Moreover, some concluding comments are also made.

\subsubsection{Computational performance}
\noindent Regarding computational performance, which can be assessed by the time an instance takes to solve and the number of iterations that were stopped due to time limit conditions being reached, the following remarks can be made:
\begin{itemize}[label={--}]
    \item Increases in the size of the instance, for the considered instance types $p(n_i.n_j.n_t.n_k.n_\ell. n_p.n_q)$, invariably lead to an increase in computational complexity;
    \item For the considered ranges, a decrease in the number of fixed roles in $e_{ijt}$ or in the unavailability percentages for the availability parameters, $l_{ik\ell}$ and $m_{k\ell p}$, also lead to notable increases in the CPU times;
    \item Conversely, the variation on the percentages for each big-$M$ upper bound, $v_{i\ell}$ and $h_{i\ell}$, did not produce such a unidimensional variation in computational complexity, with some instance types being solved more efficiently when the smaller upper bound is more frequent, and other instances seemingly showing the opposite trend.
\end{itemize}

For the first two points, the mentioned parameter variations increase the number of possible feasible variable combinations. Thus, making the instances more challenging to solve. For the irregular behaviour presented in the last remark, the explanation might be that, occasionally, the increase in the number of committee members that are assigned the parameter values $v_{i\ell}=[2,1]$ and $h_{i\ell}=[2,1]$, instead of $v_{i\ell}=[1]$ and $h_{i\ell}=[1]$, is potentially reducing solution symmetries. This effect surpasses the impact of the weaker linear relaxations induced by  larger big-$M$ upper bounds. 

\subsubsection{Number of schedulable defences}

\noindent Regarding the number of schedulable defences, logically, increasing the number of committee members and defences leads to more defences being scheduled. Additionally, the following remarks on the percentage of schedulable defences can be made:

\begin{itemize}[label={--}]

    \item Increasing the number of available rooms improves the number of schedulable defences;
    \item Reducing the number of fixed roles in $e_{ijt}$ or the unavailability percentages for the availability parameters, $l_{ik\ell}$ and $m_{k\ell p}$, promoted higher schedulability percentages.
\end{itemize}

These points show that these parameters, which are involved in the model, affect the probability that there is a suitable time slot for the scheduling of a given defence. Moreover, an increase in this probability improves the number of schedulable defences. Let us note that the same parameters that positively affect the number of schedulable defences also negatively impact the time it takes to solve the respective instance. This is explained by the increase in the number of possible assignment combinations promoted by the variation of these parameters. Conversely, the big-$M$ upper bound distribution variation does not impact the number of schedulable defences, which might indicate why its effect on the computational efficiency is not as streamlined as it is for the other considered parameters.

\subsubsection{Type of iteration distribution}
\noindent Regarding the type of iteration distribution, while it is harder to take conclusions considering the more computationally complex instances, which had iterations stop due to the set CPU time limit being reached, the following remarks can be made for the instances where these time-related stop conditions were not met:

\begin{itemize}[label={--}]
    \item The distribution of feasible and infeasible iterations,  $|N|+skip^N$ and $|I|+skip^I$, respectively, shows a slight variation, with most instances having between 75 and 85 feasible iterations out of 100;
    \item Instances that had more effective iterations took longer to solve when compared to similarly-sized instances, which had more skipped iterations. Nonetheless, their time \textit{per} effective iteration is not considerably different;
    \item The number of different non-dominated solutions shows a slight positive correlation with the percentage of schedulable defences. 
\end{itemize}

When looking at the results of the computational experiments,  some outlier instances occasionally pop out, which take longer to solve when compared to the instances that are most similar to them. Nonetheless, these usually occur due to, for some reason, the outlier instances having a relatively high number of different non-dominated solutions. This leads to fewer skipped iterations and, thus, longer CPU times. Nonetheless, this does not mean that each iteration is harder to solve, just that there are more effective iterations. Moreover, besides these occasional outlier instances, the apparent rule is that there are more different solutions in instances with a higher schedulability percentage. Accordingly, when there are more assignments, more committee members are involved. This seems to promote more possible trade-offs between the different considered objectives.

\subsubsection{Concluding comments}

\noindent The proposed method showed a remarkable capacity for finding the number of schedulable defences, $g$, always reaching optimality in the first stage. This is a helpful step for real-world problems where the decision-makers are not \textit{a priori} certain that all defences are schedulable. Furthermore, finding this parameter means that at least one feasible solution is always found.

Moreover, for almost every instance with two fixed roles, the method could map the desired subset of the Pareto front without reaching any of the defined time limits. Conversely, the same cannot be said for instances with a single fixed role, specifically for the medium and large instances, which  had several iterations being stopped due to time limits. Still, we must note that, even for these instances, the method still returns several feasible solutions and some non-dominated different ones. Thus, if applied to larger real-world instances, several different options would still be presented to the decision-maker, even if their optimality could not be proven.

Furthermore, for each instance type, the number of fixed roles was the parameter whose variation had the highest impact on both the computational performance and the number of schedulable defences. While there is no certainty that this remark will hold for different parameter ranges, it can be seen as a good indicator of the effect this scheduling decision can have on the best thesis defence scheduling method to be employed by universities with different policies. For organisations where the scheduling process occurs synchronously with the assignment of the committees, hence having less fixed roles, it might be easier to schedule all the defences, but using deterministic methods might not be suitable due to the increased computational complexity. Conversely, for organisations where the scheduling process occurs separately from the assignment of the committees, hence having more fixed roles, finding an available time slot for each defence might be more challenging. Still, it is considerably easier for a deterministic solution to their instance to be found within a reasonable time frame.
\section{Conclusions}
\noindent In this paper, we propose a MOMILP model for the single defence assignment class of the thesis defence scheduling problems. This problem consists of assigning a committee, time slot, and room to each thesis defence. Moreover, it is subject to a set of constraints that define the feasible region. Each of these constraints can be placed in one of three groups of constraints, specifically, scheduling complete committees, committee composition, or committee member and room availability. Moreover, the quality of the schedules can be assessed under two points of view, rendered operational through several criteria, specifically, committee assignment and schedule quality. To tackle the multi-objective nature of the problem, we implement an adaptation of the augmented $\epsilon-$constraint method, which allows for the mapping of a subset of the Pareto front, presenting the decision-makers with a variable number of different non-dominated solutions, while employing iteration skipping mechanisms to improve the overall computational efficiency.

The thesis defence scheduling literature has primarily been focused on solving the problem at the authors' universities. Conversely, one of the main contributions of our work is that we formalise the problem in a manner that 
is easier to adapt to institutions with different policies and regulations. Furthermore, besides including and offering a new take on the formulation of the most common objectives in thesis defence scheduling, we also regard some additional ones that were not  previously considered. Moreover, we also account for the possibility that, in instances where not all defences are schedulable, it can still be valuable for the decision-makers to be presented with an ``incomplete'' schedule so that they have access to more information and may better assess how to proceed in solving the problem.

With this work, we aim to promote the study of a fundamental academic scheduling problem, which  is remarkably underrepresented in the literature for how impactful and time-consuming it can be. Thus, we also present a novel random instance generator that can help to provide instances for future research.

The computational experiments proved that the first stage of the proposed method introduces an important step for solving thesis defence scheduling instances where it is not known that every defence is schedulable. Moreover, even for larger instances, the method was always capable of returning several solutions. Furthermore, for smaller instances or  those with two fixed roles, the optimality of each returned solution is practically always proven.

While we attempt to include most concerns and policies addressed in the literature and inclusively consider new ones, it is entirely possible and expected that some additional regulations and preferences not yet discussed but present in other universities might not have been covered by our work. Moreover, unlike other academic scheduling problems, the development and improvement of novel solution methods and algorithms are lacking for the thesis scheduling problem. Thus, we believe this to be a promising new field for future research, and that new findings can help not only the optimisation of thesis scheduling in universities but also apply to other scheduling problems, such as course timetabling or exam scheduling.
\section*{Acknowledgements}
\addcontentsline{toc}{section}{\numberline{}Acknowledgements}

\noindent João Almeida acknowledges the help of  Inês Marques, who, along with Daniel Rebelo dos Santos, introduced him to the topic of  master's thesis defence scheduling problems.  José Rui Figueira acknowledges the support by national funds through FCT (Fundação para a Ciência e a Tecnologia), under DOME Project (reference: PTDC/CCI-COM/31198/2017). All the authors acknowledge the Portuguese national funds through the FCT - Foundation for Science and Technology, I.P., under the project UIDB/00097/2020.

\addcontentsline{toc}{section}{\numberline{}References}
\bibliographystyle{model2-names}
\bibliography{references} 

\begin{thebibliography}{25}
\expandafter\ifx\csname natexlab\endcsname\relax\def\natexlab#1{#1}\fi
\expandafter\ifx\csname url\endcsname\relax
  \def\url#1{\texttt{#1}}\fi
\expandafter\ifx\csname urlprefix\endcsname\relax\def\urlprefix{URL }\fi
\providecommand{\eprint}[2][]{\url{#2}}
\providecommand{\bibinfo}[2]{#2}
\ifx\xfnm\relax \def\xfnm[#1]{\unskip,\space#1}\fi
\bibitem[{Ahmadi et~al.(2016)Ahmadi, Zandieh, Farrokh and
  Emami}]{AhmadiEtAl2016}
\bibinfo{author}{Ahmadi, E.}, \bibinfo{author}{Zandieh, M.},
  \bibinfo{author}{Farrokh, M.}, \bibinfo{author}{Emami, {\relax S.M.}.},
  \bibinfo{year}{2016}.
\newblock \bibinfo{title}{A multi objective optimization approach for flexible
  job shop scheduling problem under random machine breakdown by evolutionary
  algorithms}.
\newblock \bibinfo{journal}{Computers and Operations Research}
  \bibinfo{volume}{73}, \bibinfo{pages}{56--66}.
\bibitem[{Al-Yakoob and Sherali(2015)}]{Al-YakoobSherali2015}
\bibinfo{author}{Al-Yakoob, {\relax S.M.}.}, \bibinfo{author}{Sherali, {\relax
  H.D.}.}, \bibinfo{year}{2015}.
\newblock \bibinfo{title}{Mathematical models and algorithms for a high school
  timetabling problem}.
\newblock \bibinfo{journal}{Computers and Operations Research}
  \bibinfo{volume}{61}, \bibinfo{pages}{56--68}.
\bibitem[{Babaei et~al.(2015)Babaei, Karimpour and Hadidi}]{BabaeiEtAl2015}
\bibinfo{author}{Babaei, H.}, \bibinfo{author}{Karimpour, J.},
  \bibinfo{author}{Hadidi, A.}, \bibinfo{year}{2015}.
\newblock \bibinfo{title}{A survey of approaches for university course
  timetabling problem}.
\newblock \bibinfo{journal}{Computers and Industrial Engineering}
  \bibinfo{volume}{86}, \bibinfo{pages}{43--59}.
\bibitem[{Battistutta et~al.(2019)Battistutta, Ceschia, De~Cesco, Di~Gaspero
  and Schaerf}]{BattistuttaEtAl2019}
\bibinfo{author}{Battistutta, M.}, \bibinfo{author}{Ceschia, S.},
  \bibinfo{author}{De~Cesco, F.}, \bibinfo{author}{Di~Gaspero, L.},
  \bibinfo{author}{Schaerf, A.}, \bibinfo{year}{2019}.
\newblock \bibinfo{title}{Modelling and solving the thesis defense timetabling
  problem}.
\newblock \bibinfo{journal}{Journal of the Operational Research Society}
  \bibinfo{volume}{70}, \bibinfo{pages}{1039--1050}.
\bibitem[{Bellio et~al.(2021)Bellio, Ceschia, Di~Gaspero and
  Schaerf}]{BellioEtAl2021}
\bibinfo{author}{Bellio, R.}, \bibinfo{author}{Ceschia, S.},
  \bibinfo{author}{Di~Gaspero, L.}, \bibinfo{author}{Schaerf, A.},
  \bibinfo{year}{2021}.
\newblock \bibinfo{title}{Two-stage multi-neighborhood simulated annealing for
  uncapacitated examination timetabling}.
\newblock \bibinfo{journal}{Computers and Operations Research}
  \bibinfo{volume}{132}.
\bibitem[{Borgonjon and Maenhout(2022)}]{BorgonjonMaenhout2022}
\bibinfo{author}{Borgonjon, T.}, \bibinfo{author}{Maenhout, B.},
  \bibinfo{year}{2022}.
\newblock \bibinfo{title}{A two-phase pareto front method for solving the
  bi-objective personnel task rescheduling problem}.
\newblock \bibinfo{journal}{Computers and Operations Research}
  \bibinfo{volume}{138}.
\bibitem[{Burke et~al.(2010)Burke, Mareček, Parkes and
  Rudová}]{BurkeEtAl2010Decomposition}
\bibinfo{author}{Burke, {\relax E.K}.}, \bibinfo{author}{Mareček, J.},
  \bibinfo{author}{Parkes, {\relax A.J}.}, \bibinfo{author}{Rudová, H.},
  \bibinfo{year}{2010}.
\newblock \bibinfo{title}{Decomposition, reformulation, and diving in
  university course timetabling}.
\newblock \bibinfo{journal}{Computers and Operations Research}
  \bibinfo{volume}{37}, \bibinfo{pages}{582--597}.
\newblock \bibinfo{note}{Hybrid Metaheuristics}.
\bibitem[{Chaudhuri and De(2010)}]{ChaudhuriDe2010}
\bibinfo{author}{Chaudhuri, A.}, \bibinfo{author}{De, K.},
  \bibinfo{year}{2010}.
\newblock \bibinfo{title}{Fuzzy genetic heuristic for university course
  timetable problem}.
\newblock \bibinfo{journal}{International Journal of Advances in Soft Computing
  and its Applications} \bibinfo{volume}{2}, \bibinfo{pages}{100--123}.
\bibitem[{Christopher and Wicaksana(2021)}]{ChristopherWicaksana2021}
\bibinfo{author}{Christopher, G.}, \bibinfo{author}{Wicaksana, A.},
  \bibinfo{year}{2021}.
\newblock \bibinfo{title}{Particle swarm optimization for solving thesis
  defense timetabling problem}.
\newblock \bibinfo{journal}{TELKOMNIKA Telecommunication, Computing,
  Electronics and Control} \bibinfo{volume}{19}, \bibinfo{pages}{762--769}.
\bibitem[{Goh et~al.(2017)Goh, Kendall and Sabar}]{GohEtAl2017}
\bibinfo{author}{Goh, {\relax S.L.}.}, \bibinfo{author}{Kendall, G.},
  \bibinfo{author}{Sabar, {\relax N.R.}.}, \bibinfo{year}{2017}.
\newblock \bibinfo{title}{Improved local search approaches to solve the post
  enrolment course timetabling problem}.
\newblock \bibinfo{journal}{European Journal of Operational Research}
  \bibinfo{volume}{261}, \bibinfo{pages}{17--29}.
\bibitem[{Huynh et~al.(2012)Huynh, Pham and Pham}]{HuynhEtAl2012}
\bibinfo{author}{Huynh, {\relax T. T. B}.}, \bibinfo{author}{Pham, {\relax Q.
  D}.}, \bibinfo{author}{Pham, {\relax D. D}.}, \bibinfo{year}{2012}.
\newblock \bibinfo{title}{Genetic algorithm for solving the master thesis
  timetabling problem with multiple objectives}, in: \bibinfo{booktitle}{TAAI
  '12: Proceedings of the 2012 Conference on Technologies and Applications of
  Artificial Intelligence}, \bibinfo{publisher}{IEEE Computer Society},
  \bibinfo{address}{Washington, D.C., USA}. pp. \bibinfo{pages}{74--79}.
\bibitem[{Kochaniková and Rudová(2013)}]{KochanikovRudova2013}
\bibinfo{author}{Kochaniková, B.}, \bibinfo{author}{Rudová, H.},
  \bibinfo{year}{2013}.
\newblock \bibinfo{title}{Student scheduling for bachelor state examinations},
  in: \bibinfo{editor}{Kendall, G.}, \bibinfo{editor}{Berghe, G.V.},
  \bibinfo{editor}{McCollum, B.} (Eds.), \bibinfo{booktitle}{6th
  Multidisciplinary International Conference on Scheduling : Theory and
  Applications}, \bibinfo{publisher}{Springer}, \bibinfo{address}{New York, NY,
  USA}.
\bibitem[{Lei et~al.(2016)Lei, Wang, Zhang, Liu and Zha}]{LeiEtAl2016}
\bibinfo{author}{Lei, H.}, \bibinfo{author}{Wang, R.}, \bibinfo{author}{Zhang,
  T.}, \bibinfo{author}{Liu, Y.}, \bibinfo{author}{Zha, Y.},
  \bibinfo{year}{2016}.
\newblock \bibinfo{title}{A multi-objective co-evolutionary algorithm for
  energy-efficient scheduling on a green data center}.
\newblock \bibinfo{journal}{Computers and Operations Research}
  \bibinfo{volume}{75}, \bibinfo{pages}{103--117}.
\bibitem[{Limanto et~al.(2019)Limanto, Benarkah and Adelia}]{LimantoEtAl2019}
\bibinfo{author}{Limanto, S.}, \bibinfo{author}{Benarkah, N.},
  \bibinfo{author}{Adelia, T.}, \bibinfo{year}{2019}.
\newblock \bibinfo{title}{Thesis examination timetabling using genetic
  algorithm}, in: \bibinfo{editor}{Muliawati, T.}, \bibinfo{editor}{Ardiansyah,
  M.}, \bibinfo{editor}{Sari, D.}, \bibinfo{editor}{Permatasari, D.},
  \bibinfo{editor}{Mu'arifin} (Eds.), \bibinfo{booktitle}{International
  Electronics Symposium on Knowledge Creation and Intelligent Computing,
  IES-KCIC 2018 - Proceedings}, \bibinfo{publisher}{Institute of Electrical and
  Electronics Engineers Inc.}, \bibinfo{address}{Piscataway, NJ, USA}. pp.
  \bibinfo{pages}{6--10}.
\bibitem[{Liu et~al.(2014)Liu, Zhang, Zhu and Rao}]{LiuEtAl2014}
\bibinfo{author}{Liu, Q.}, \bibinfo{author}{Zhang, C.}, \bibinfo{author}{Zhu,
  K.}, \bibinfo{author}{Rao, Y.}, \bibinfo{year}{2014}.
\newblock \bibinfo{title}{Novel multi-objective resource allocation and
  activity scheduling for fourth party logistics}.
\newblock \bibinfo{journal}{Computers and Operations Research}
  \bibinfo{volume}{44}, \bibinfo{pages}{42--51}.
\bibitem[{Mavrotas(2009)}]{Mavrotas2009}
\bibinfo{author}{Mavrotas, G.}, \bibinfo{year}{2009}.
\newblock \bibinfo{title}{Effective implementation of the e-constraint method
  in multi-objective mathematical programming problems}.
\newblock \bibinfo{journal}{Applied Mathematics and Computation}
  \bibinfo{volume}{213}, \bibinfo{pages}{455--465}.
\bibitem[{Mavrotas and Florios(2013)}]{MavrotasFlorios2013}
\bibinfo{author}{Mavrotas, G.}, \bibinfo{author}{Florios, K.},
  \bibinfo{year}{2013}.
\newblock \bibinfo{title}{An improved version of the augmented e-constraint
  method {(AUGMECON2)} for finding the exact pareto set in multi-objective
  integer programming problems}.
\newblock \bibinfo{journal}{Applied Mathematics and Computation}
  \bibinfo{volume}{219}, \bibinfo{pages}{9652--9669}.
\bibitem[{Pham et~al.(2015)Pham, Hoang, Huynh and Nguyen}]{PhamEtAl2015}
\bibinfo{author}{Pham, {\relax Q.D}.}, \bibinfo{author}{Hoang, {\relax
  T.A.D}.}, \bibinfo{author}{Huynh, {\relax T.T}.}, \bibinfo{author}{Nguyen,
  {\relax T.H}.}, \bibinfo{year}{2015}.
\newblock \bibinfo{title}{A {Java} library for constraint-based local search:
  Application to the master thesis defense timetabling problem}, in:
  \bibinfo{editor}{Quyet~Tang, H.}, \bibinfo{editor}{Le~Anh, P.},
  \bibinfo{editor}{De~Raedt, L.}, \bibinfo{editor}{Deville, Y.},
  \bibinfo{editor}{Bui, M.}, \bibinfo{editor}{Trong Thi~Dieu, L.},
  \bibinfo{editor}{Nguyen~Thi, O.}, \bibinfo{editor}{Dinh~Viet, S.},
  \bibinfo{editor}{Nguyen~Ba, N.} (Eds.), \bibinfo{booktitle}{Proceedings of
  the Sixth International Symposium on Information and Communication
  Technology}, \bibinfo{publisher}{Association for Computing Machinery},
  \bibinfo{address}{New York, NY, USA}. pp. \bibinfo{pages}{67--74}.
\bibitem[{Santiago-Mozos et~al.(2005)Santiago-Mozos, Salcedo-Sanz,
  Deprado-Cumplido and Bousoño-Calzón}]{Santiago-MozosEtAl2005}
\bibinfo{author}{Santiago-Mozos, R.}, \bibinfo{author}{Salcedo-Sanz, S.},
  \bibinfo{author}{Deprado-Cumplido, M.}, \bibinfo{author}{Bousoño-Calzón,
  C.}, \bibinfo{year}{2005}.
\newblock \bibinfo{title}{A two-phase heuristic evolutionary algorithm for
  personalizing course timetables: A case study in a spanish university}.
\newblock \bibinfo{journal}{Computers and Operations Research}
  \bibinfo{volume}{32}, \bibinfo{pages}{1761--1776}.
\bibitem[{Su et~al.(2020)Su, Luo, Deng, Xia and Guo}]{SuEtAl2020}
\bibinfo{author}{Su, P.}, \bibinfo{author}{Luo, {\relax B.X}.},
  \bibinfo{author}{Deng, {\relax F.Y}.}, \bibinfo{author}{Xia, {\relax A.X}.},
  \bibinfo{author}{Guo, Y.}, \bibinfo{year}{2020}.
\newblock \bibinfo{title}{Group strategy of dissertation defense based on
  greedy retrospective hybrid algorithm}, in: \bibinfo{booktitle}{Journal of
  Physics: Conference Series}, \bibinfo{publisher}{IOP Publishing Ltd},
  \bibinfo{address}{Bristol, UK}.
\bibitem[{Sørensen and Dahms(2014)}]{SorensenDahms2014}
\bibinfo{author}{Sørensen, M.}, \bibinfo{author}{Dahms, F.},
  \bibinfo{year}{2014}.
\newblock \bibinfo{title}{A two-stage decomposition of high school timetabling
  applied to cases in denmark}.
\newblock \bibinfo{journal}{Computers and Operations Research}
  \bibinfo{volume}{43}, \bibinfo{pages}{36--49}.
\bibitem[{Tan et~al.(2021)Tan, Goh, Kendall and Sabar}]{TanEtAl2021}
\bibinfo{author}{Tan, {\relax J.S.}.}, \bibinfo{author}{Goh, {\relax S.L.}.},
  \bibinfo{author}{Kendall, G.}, \bibinfo{author}{Sabar, {\relax N.R.}.},
  \bibinfo{year}{2021}.
\newblock \bibinfo{title}{A survey of the state-of-the-art of optimisation
  methodologies in school timetabling problems}.
\newblock \bibinfo{journal}{Expert Systems with Applications}
  \bibinfo{volume}{165}.
\bibitem[{Tawakkal and Suyanto(2020)}]{TawakkalSuyanto2020}
\bibinfo{author}{Tawakkal, {\relax M.I}.}, \bibinfo{author}{Suyanto},
  \bibinfo{year}{2020}.
\newblock \bibinfo{title}{Exploration-exploitation balanced krill herd
  algorithm for thesis examination timetabling}, in: \bibinfo{booktitle}{2020
  International Conference on Data Science and Its Applications},
  \bibinfo{publisher}{Institute of Electrical and Electronics Engineers Inc.},
  \bibinfo{address}{Piscataway, NJ, USA}.
\bibitem[{Vermuyten et~al.(2016)Vermuyten, Lemmens, Marques and
  Beliën}]{VermuytenEtAl2016}
\bibinfo{author}{Vermuyten, H.}, \bibinfo{author}{Lemmens, S.},
  \bibinfo{author}{Marques, I.}, \bibinfo{author}{Beliën, J.},
  \bibinfo{year}{2016}.
\newblock \bibinfo{title}{Developing compact course timetables with optimized
  student flows}.
\newblock \bibinfo{journal}{European Journal of Operational Research}
  \bibinfo{volume}{251}, \bibinfo{pages}{651--661}.
\bibitem[{Wang et~al.(2017)Wang, Liu and Jin}]{WangEtAl2017}
\bibinfo{author}{Wang, {\relax D.-J.}.}, \bibinfo{author}{Liu, F.},
  \bibinfo{author}{Jin, Y.}, \bibinfo{year}{2017}.
\newblock \bibinfo{title}{A multi-objective evolutionary algorithm guided by
  directed search for dynamic scheduling}.
\newblock \bibinfo{journal}{Computers and Operations Research}
  \bibinfo{volume}{79}, \bibinfo{pages}{279--290}.

\end{thebibliography}

\vfill\newpage

\appendix

\section{Mathematical estimation of the availability and unavailability percentages}

\noindent This appendix refers to an alternative way of estimating the probability distribution for the different considered availability states, $\alpha=0,\ldots,n_\alpha$. It became redundant when we made the connection to Markov Chains and, for this reason, it is present here only as a curiosity.

Let us note that in this appendix we differentiate the conditional probabilities that are an input, $p(\alpha|\alpha)^{inp}$, and those that are the actual conditional probabilities, $p(\alpha|\alpha)$. This differentiation is necessary due to the exceptional additions.

The expected duration of the availability blocks for a $l_{ik\ell} = \alpha$ or a $m_{k\ell p} = \alpha$, $d^a_\alpha, \; \alpha\neq0$ can be computed through Equation \eqref{available duration}, and that of the unavailability blocks, with $l_{ik\ell} = 0$ or $m_{k\ell p} = 0$, $d_0$, through Equation \eqref{unavailable duration}.

\begin{equation}\label{available duration}
\displaystyle
    E(d^a_\alpha) =(1 - p(\alpha|\alpha)^{inp})^{-1},\;\alpha\neq0
\end{equation}

\begin{equation}\label{unavailable duration}
\displaystyle
    E(d_0) = d - 1 + (1 - p(0|0)^{inp})^{-1}
\end{equation}

\begin{theorem}
The expected duration of an availability block, for a given $l_{ik\ell} = \alpha$ or $m_{k\ell p} = \alpha$, $d^a_\alpha$, is determined by Equation \eqref{available duration}.
\end{theorem}
\begin{proof}

An availability block is generated through consecutive Bernoulli trials, with a probability of a successful change $p=1 - p(\alpha|\alpha)^{inp}$, and a non-success from a previous set of trials. Let $d^B_\alpha$ denote the number of unsuccessful Bernoulli trials before a success occurs. Given that the number of Bernoulli trials until a success occurs follows a geometric distribution, $G(1 - p(\alpha|\alpha)^{inp})$, the expected $d^B_\alpha$ is computed as follows,

\begin{center}
    $E(d^B_\alpha) = (1 - p(\alpha|\alpha)^{inp})^{-1} - 1$.
\end{center}

Since each availability block contains both its correspondent non-successes, and the previous success, then

\begin{center}
    $E(d^a_\alpha) = 1 + E(d^B_\alpha)$.
\end{center}

Which proves the result of Equation \eqref{available duration}.
\end{proof}

\begin{corollary}
The average duration of an unavailability block is determined by Equation \eqref{unavailable duration}. 
\begin{proof}
On the one hand, the number, $d-1$, of unavailable time-slots is added \textit{per} unavailability block. On the other hand, the remaining generation process is similar to that of the availability blocks, hence, the same principles from Theorem 1 apply, proving the result of Equation \eqref{unavailable duration}.
\end{proof}
\end{corollary}

To conduct the computational experiments, it is also of interest to be able to estimate the probability $p(\alpha)$ that a given time-slot will be of a given type $\alpha=0,\ldots,n_\alpha$, given the input conditional probabilities, $p(\alpha|\alpha)^{inp}$. Moreover, we can assume that, given a set of inputs, this distribution tends to a uniform set of values for $p(\alpha),\; \alpha=0,\ldots, n_\alpha$, that is, for any $\ell$, the probability of assigning the value $\alpha$ is $p({\alpha^\ell}) = p({\alpha^{\overline{\ell}}}) = p(\alpha)$. Knowing that the sum of all probabilities $p(\alpha)$ must be equal to 1, the distribution can be found by solving System \eqref{equation system}.

\begin{equation}\label{equation system}
    \begin{cases}
    \displaystyle
      p({0}) = \sum_{\alpha=0}^{n_\alpha} p({\alpha})p(0|\alpha)\\
      \hspace{1.4cm}\vdots\\
      \displaystyle
      p({\overline{\alpha}}) = \sum_{\alpha=0}^{n_\alpha} p(\alpha)p(\overline{\alpha}|\alpha)\\
      \hspace{1.4cm}\vdots\\
      \displaystyle
      p({n_\alpha}) = \sum_{\alpha=0}^{n_\alpha} p(\alpha)p(n_\alpha|\alpha)\\
      \displaystyle
      \sum_{\alpha=0}^{n_\alpha} p(\alpha) = 1\\
    \end{cases}
\end{equation}

To compute the set of probabilities $p(\alpha)$ through System \eqref{equation system}, we must know all conditional probabilities, $p(\alpha|\overline{\alpha})$:

\begin{lemma}
 All conditional probabilities, $p(\alpha|\overline{\alpha})$, are defined, given a certain set of input conditional probabilities, $p(\alpha|\alpha)^{inp}$.
\end{lemma}

\begin{proof}
The values for the success probabilities $p(\alpha|\overline{\alpha})^{inp},\; \overline{\alpha}\neq0, \alpha\neq\overline{\alpha}$, can be computed, based on the inputs, through Equation \eqref{equation conditional change}. However, when the exception conditions are met, $p(0|0_e) = 1$, and $p(\alpha|0_e) = 0,\;\alpha\neq0$. Unless there are no exceptional additions, which is only true if $d-1=0$, it follows that $p(\alpha|0)\neq p(\alpha|0)^{inp}$. 

Let $S$ denote the sample space of all the generated unavailability blocks, $\alpha = 0$. An $\alpha\in S$ will either lead to an exceptional addition, with probability $p(e)$, or, to a Bernoulli trial, with a probability $1-p(e)$ and a success probability of $1-p(0|0)^{inp}$. Hence, these are pairwise disjoint events, whose union will define the sample space, $S$.

Thus, applying the law of total probability and the conditional probability definition,

\begin{center}
    $p(\alpha|0) = p(e)p(\alpha|0_e) + (1-p(e))p(\alpha|0)^{inp}$.
\end{center}

Now we need to find the probability, $p(e)$. Following Corollary 1, we know the average duration, $d_0$, of an unavailability block. Thus, considering $n_b$ blocks, the total number of $\alpha=0$, $n_0$ can be defined as follows,

\begin{center}
    $n_0 = d_0n_b$.
\end{center}

Moreover, we also know that each block invariably contains $d-1$ exceptionally added $\alpha=0$, $0_e$, thus, the total number of $0_e$ that will lead to an exceptional addition, $n_{0_e}$ can be defined as follows,

\begin{center}
    $n_{0_e} = (d-1)n_b$.
\end{center}

Thus, by definition, $p(e)$ is the ratio between $n_{0_e}$ and $n_0$, hence,

\begin{center}
    $p(e) = (d-1)(d_0)^{-1}$.
\end{center}

Therefore, we know all parameters that are necessary to compute all conditional probabilities, $p(\alpha|\overline{\alpha})$.
\end{proof}

Now, to prove the feasibility of System \eqref{equation system}, i.e., the system represents the probabilities $p(0),\ldots,p({n_\alpha})$, of any parameter $l_{ik\ell}$ or $m_{k\ell p}$ being assigned the values $\alpha=0,\ldots,\alpha=n_\alpha$:

\begin{theorem}
Assuming that, $p({\alpha^\ell}) = p({\alpha^{\overline{\ell}}}) = p(\alpha)$, the probability distribution for $p(\alpha)$ can be computed through system \eqref{equation system}.
\end{theorem}

\begin{proof}
Let $\hat{S}$ denote the sample space of randomly generated availability parameters, $\hat{p}$, each of them considering time-slots, $\ell \in \{1,\ldots,n_\ell\}$. Let $\ell$ and $\overline{\ell} = \ell - 1$, denote two consecutive hour-slots, respectively. Let $\overline{S}$ denote the event space of possible events, $\alpha^\ell \in \{0,\ldots,n_\alpha\}$, for any given time-slot $\ell$. By definition, $\overline{S}$ is composed of $n_\alpha + 1$ pairwise disjoint events of probability $p(\alpha)$. Moreover, their union includes the entire sample space $\hat{S}$, because, for every $\hat{p}\in \hat{S}$ and $\ell$, one $\alpha^\ell \in \overline{S}$ is assigned. Hence, we can apply the law of total probability, and represent the probability $p({\alpha^\ell})$ as follows,

\begin{equation*}
    {\displaystyle
    p({\alpha^\ell}) = \sum_{\overline\alpha}=0}^{n_\alpha} p({\alpha^\ell})\cap p({\overline{\alpha}^{\overline{\ell}}}),\;\alpha = 0, \ldots, n_\alpha, \overline{\ell} = \ell-1
\end{equation*}

Applying the definition of conditional probability and holding our assumption that $p({\alpha^\ell}) = p({\alpha^{\overline{\ell}}}) = p(\alpha)$ as true, the probability $p(\alpha)$ of a parameter value $\alpha$ occurring at any time-slot can be defined as follows,

\begin{equation*}\label{substitution}
\displaystyle
    p(\alpha) = \sum_{\overline{\alpha}=0}^{n_\alpha} p({\overline{\alpha}})p(\alpha|\overline{\alpha}),\;\alpha = 0, \ldots, n_\alpha
\end{equation*}

Therefore, and knowing that the sum of all $p(\alpha)$ must be equal to $1$, the feasibility of system \eqref{equation system} is proven. Moreover, through Lemma 1, we know that all $p(\alpha|\overline{\alpha})$ are known, hence System \eqref{equation system} can be solved to find the distribution of probabilities, $p(\alpha)$.
\end{proof}

\newpage
\section{Computational results}

\begin{landscape}
\begin{table}[]
    \caption{Computational experiments - small instances ($n_i = 25, \; n_j = 20$)}\label{small}
    \scalebox{0.8}{
    \centering
    \begin{tabular}{|c|c|cccccccccc|cccccccc|}
    \hline
    \multicolumn{2}{|c|}{{Instance}}&\multicolumn{10}{c|}{{Data}}&\multicolumn{8}{c|}{{Output}}\\
    \hline
    N&$p(n_i.n_j.n_t.n_k.n_\ell. n_p.n_q)$&$d$&$u_i$&$e_{ijt}$&$c_i$&$l_{ik\ell}$&$m_{k\ell p}$&$v_{i\ell}$&$h_{i\ell}$&$r_{iq}$&$t_{iq}$&$|N|$&$|I|$&$skip^N$&$skip^I$&$time^N$&$time^I$&$g$&CPU\\
    \hline
1&$p(25.20.3.15.16.3.15)$ & 2 & [0.7, 0.3] & 2 & 13& [0.82, 0.09, 0.09] & [0.86, 0.14] & [0.8, 0.2] &[0.8, 0.2] & 3 & 3 &21&3&65&11&0&0&10&230\\
2&$p(25.20.3.15.16.3.15)$ & 2 & [0.7, 0.3] & 2 & 13& [0.82, 0.09, 0.09] & [0.86, 0.14] & [0.7, 0.3] &[0.7, 0.3] & 3 & 3 &26&5&53&16&0&0&18&289\\
3&$p(25.20.3.15.16.3.15)$ & 2 & [0.7, 0.3] & 2 & 13& [0.82, 0.09, 0.09] & [0.8, 0.2] & [0.8, 0.2] &[0.8, 0.2] & 3 & 3 &25&4&64&7&0&0&17&550\\
4&$p(25.20.3.15.16.3.15)$ & 2 & [0.7, 0.3] & 2 & 13& [0.82, 0.09, 0.09] & [0.8, 0.2] & [0.7, 0.3] &[0.7, 0.3] & 3 & 3 &11&1&86&2&0&0&14&267\\
5&$p(25.20.3.15.16.3.15)$ & 2 & [0.7, 0.3] & 2 & 13& [0.78, 0.11, 0.11] & [0.86, 0.14] & [0.8, 0.2] &[0.8, 0.2] & 3 & 3 &48&5&30&17&0&0&20&542\\
6&$p(25.20.3.15.16.3.15)$ & 2 & [0.7, 0.3] & 2 & 13& [0.78, 0.11, 0.11] & [0.86, 0.14] & [0.7, 0.3] &[0.7, 0.3] & 3 & 3 &51&5&27&17&0&0&18&1067\\
7&$p(25.20.3.15.16.3.15)$ & 2 & [0.7, 0.3] & 2 & 13& [0.78, 0.11, 0.11] & [0.8, 0.2] & [0.8, 0.2] &[0.8, 0.2] & 3 & 3 &59&5&22&14&0&0&17&590\\
8&$p(25.20.3.15.16.3.15)$ & 2 & [0.7, 0.3] & 2 & 13& [0.78, 0.11, 0.11] & [0.8, 0.2] & [0.7, 0.3] &[0.7, 0.3] & 3 & 3 &38&5&44&13&0&0&19&1135\\
9&$p(25.20.3.15.16.3.15)$ & 2 & [0.7, 0.3] & 1 & 13& [0.86, 0.07, 0.07] & [0.86, 0.14] & [0.8, 0.2] &[0.8, 0.2] & 3 & 3 &47&5&32&16&0&0&18&513\\
10&$p(25.20.3.15.16.3.15)$ & 2 & [0.7, 0.3] & 1 & 13& [0.86, 0.07, 0.07] & [0.86, 0.14] & [0.7, 0.3] &[0.7, 0.3] & 3 & 3 &25&6&49&20&0&0&14&298\\
11&$p(25.20.3.15.16.3.15)$ & 2 & [0.7, 0.3] & 1 & 13& [0.86, 0.07, 0.07] & [0.8, 0.2] & [0.8, 0.2] &[0.8, 0.2] & 3 & 3 &44&5&33&18&0&0&20&1418\\
12&$p(25.20.3.15.16.3.15)$ & 2 & [0.7, 0.3] & 1 & 13& [0.86, 0.07, 0.07] & [0.8, 0.2] & [0.7, 0.3] &[0.7, 0.3] & 3 & 3 &38&5&36&21&0&0&20&417\\
13&$p(25.20.3.15.16.3.15)$ & 2 & [0.7, 0.3] & 1 & 13& [0.82, 0.09, 0.09] & [0.86, 0.14] & [0.8, 0.2] &[0.8, 0.2] & 3 & 3 &42&6&34&18&0&0&20&485\\
14&$p(25.20.3.15.16.3.15)$ & 2 & [0.7, 0.3] & 1 & 13& [0.82, 0.09, 0.09] & [0.86, 0.14] & [0.7, 0.3] &[0.7, 0.3] & 3 & 3 &50&5&29&16&0&0&20&777\\
15&$p(25.20.3.15.16.3.15)$ & 2 & [0.7, 0.3] & 1 & 13& [0.82, 0.09, 0.09] & [0.8, 0.2] & [0.8, 0.2] &[0.8, 0.2] & 3 & 3 &57&7&12&24&0&0&19&1140\\
16&$p(25.20.3.15.16.3.15)$ & 2 & [0.7, 0.3] & 1 & 13& [0.82, 0.09, 0.09] & [0.8, 0.2] & [0.7, 0.3] &[0.7, 0.3] & 3 & 3 &47&6&27&20&0&0&20&10246\\

\hdashline
17&$p(25.20.3.15.16.4.15)$ & 2 & [0.7, 0.3] & 2 & 13& [0.82, 0.09, 0.09]  & [0.86, 0.14] & [0.8, 0.2] &[0.8, 0.2] & 3 & 3 &33&5&41&21&0&0&16&488\\
18&$p(25.20.3.15.16.4.15)$ & 2 & [0.7, 0.3] & 2 & 13& [0.82, 0.09, 0.09]  & [0.86, 0.14] & [0.7, 0.3] &[0.7, 0.3] & 3 & 3 &25&3&56&16&0&0&15&384\\
19&$p(25.20.3.15.16.4.15)$ & 2 & [0.7, 0.3] & 2 & 13& [0.82, 0.09, 0.09]  & [0.8, 0.2] & [0.8, 0.2] &[0.8, 0.2] & 3 & 3 &52&4&32&12&0&0&17&818\\
20&$p(25.20.3.15.16.4.15)$ & 2 & [0.7, 0.3] & 2 & 13& [0.82, 0.09, 0.09]  & [0.8, 0.2] & [0.7, 0.3] &[0.7, 0.3] & 3 & 3 &60&6&13&21&0&0&20&1053\\
21&$p(25.20.3.15.16.4.15)$ & 2 & [0.7, 0.3] & 2 & 13& [0.78, 0.11, 0.11] & [0.86, 0.14] & [0.8, 0.2] &[0.8, 0.2] & 3 & 3 &51&6&21&22&0&0&20&750\\
22&$p(25.20.3.15.16.4.15)$ & 2 & [0.7, 0.3] & 2 & 13& [0.78, 0.11, 0.11] & [0.86, 0.14] & [0.7, 0.3] &[0.7, 0.3] & 3 & 3 &24&5&55&16&0&0&13&369\\
23&$p(25.20.3.15.16.4.15)$ & 2 & [0.7, 0.3] & 2 & 13& [0.78, 0.11, 0.11] & [0.8, 0.2] & [0.8, 0.2] &[0.8, 0.2] & 3 & 3 &61&5&21&13&0&0&20&991\\
24&$p(25.20.3.15.16.4.15)$ & 2 & [0.7, 0.3] & 2 & 13& [0.78, 0.11, 0.11] & [0.8, 0.2] & [0.7, 0.3] &[0.7, 0.3] & 3 & 3 &60&4&24&12&0&0&20&2631\\
25&$p(25.20.3.15.16.4.15)$ & 2 & [0.7, 0.3] & 1 & 13& [0.86, 0.07, 0.07] & [0.86, 0.14] & [0.8, 0.2] &[0.8, 0.2] & 3 & 3 &38&6&35&21&0&0&19&1217\\
26&$p(25.20.3.15.16.4.15)$ & 2 & [0.7, 0.3] & 1 & 13& [0.86, 0.07, 0.07] & [0.86, 0.14] & [0.7, 0.3] &[0.7, 0.3] & 3 & 3 &45&5&32&18&0&0&20&1004\\
27&$p(25.20.3.15.16.4.15)$ & 2 & [0.7, 0.3] & 1 & 13& [0.86, 0.07, 0.07] & [0.8, 0.2] & [0.8, 0.2] &[0.8, 0.2] & 3 & 3 &47&5&31&17&0&0&20&4364\\
28&$p(25.20.3.15.16.4.15)$ & 2 & [0.7, 0.3] & 1 & 13& [0.86, 0.07, 0.07] & [0.8, 0.2] & [0.7, 0.3] &[0.7, 0.3] & 3 & 3 &53&6&19&22&0&0&20&1522\\
29&$p(25.20.3.15.16.4.15)$ & 2 & [0.7, 0.3] & 1 & 13& [0.82, 0.09, 0.09]  & [0.86, 0.14] & [0.8, 0.2] &[0.8, 0.2] & 3 & 3 &42&5&37&16&0&0&20&1221\\
30&$p(25.20.3.15.16.4.15)$ & 2 & [0.7, 0.3] & 1 & 13& [0.82, 0.09, 0.09]  & [0.86, 0.14] & [0.7, 0.3] &[0.7, 0.3] & 3 & 3 &38&6&37&19&0&0&20&709\\
31&$p(25.20.3.15.16.4.15)$ & 2 & [0.7, 0.3] & 1 & 13& [0.82, 0.09, 0.09]  & [0.8, 0.2] & [0.8, 0.2] &[0.8, 0.2] & 3 & 3 &40&6&32&22&0&0&20&12951\\
32&$p(25.20.3.15.16.4.15)$ & 2 & [0.7, 0.3] & 1 & 13& [0.82, 0.09, 0.09]  & [0.8, 0.2] & [0.7, 0.3] &[0.7, 0.3] & 3 & 3 &41&5&37&15&0&2&20&21175\\

    \hline
    \end{tabular}
}
    
\end{table}

\begin{table}[]
    \caption{Computational experiments - medium instances ($n_i = 38, \; n_j = 30$)}\label{medium}
    \scalebox{0.8}{
    \centering
    \begin{tabular}{|c|c|cccccccccc|cccccccc|}
    \hline
    \multicolumn{2}{|c|}{{Instance}}&\multicolumn{10}{c|}{{Data}}&\multicolumn{8}{c|}{{Output}}\\
    \hline
    N&$p(n_i.n_j.n_t.n_k.n_\ell. n_p.n_q)$&$d$&$u_i$&$e_{ijt}$&$c_i$&$l_{ik\ell}$&$m_{k\ell p}$&$v_{i\ell}$&$h_{i\ell}$&$r_{iq}$&$t_{iq}$&$|N|$&$|I|$&$skip^N$&$skip^I$&$time^N$&$time^I$&$g$&CPU\\
    \hline
33&$p(38.30.3.15.16.3.15)$ & 2 & [0.7, 0.3] & 2 & 19& [0.82, 0.09, 0.09]  & [0.86, 0.14] & [0.8, 0.2] &[0.8, 0.2] & 3 & 3 &37&6&37&20&0&0&26&973\\
34&$p(38.30.3.15.16.3.15)$ & 2 & [0.7, 0.3] & 2 & 19& [0.82, 0.09, 0.09]  & [0.86, 0.14] & [0.7, 0.3] &[0.7, 0.3] & 3 & 3 &29&2&68&1&0&0&20&698\\
35&$p(38.30.3.15.16.3.15)$ & 2 & [0.7, 0.3] & 2 & 19& [0.82, 0.09, 0.09]  & [0.8, 0.2] & [0.8, 0.2] &[0.8, 0.2] & 3 & 3 &39&5&41&15&0&0&24&1003\\
36&$p(38.30.3.15.16.3.15)$ & 2 & [0.7, 0.3] & 2 & 19& [0.82, 0.09, 0.09]  & [0.8, 0.2] & [0.7, 0.3] &[0.7, 0.3] & 3 & 3 &60&4&24&12&0&0&26&1484\\
37&$p(38.30.3.15.16.3.15)$ & 2 & [0.7, 0.3] & 2 & 19& [0.78, 0.11, 0.11] & [0.86, 0.14] & [0.8, 0.2] &[0.8, 0.2] & 3 & 3 &30&5&48&17&0&0&22&788\\
38&$p(38.30.3.15.16.3.15)$ & 2 & [0.7, 0.3] & 2 & 19& [0.78, 0.11, 0.11] & [0.86, 0.14] & [0.7, 0.3] &[0.7, 0.3] & 3 & 3 &51&5&28&16&0&0&22&1593\\
39&$p(38.30.3.15.16.3.15)$ & 2 & [0.7, 0.3] & 2 & 19& [0.78, 0.11, 0.11] & [0.8, 0.2] & [0.8, 0.2] &[0.8, 0.2] & 3 & 3 &43&4&41&12&0&0&26&1169\\
40&$p(38.30.3.15.16.3.15)$ & 2 & [0.7, 0.3] & 2 & 19& [0.78, 0.11, 0.11] & [0.8, 0.2] & [0.7, 0.3] &[0.7, 0.3] & 3 & 3 &41&5&38&16&0&0&28&10129\\
41&$p(38.30.3.15.16.3.15)$ & 2 & [0.7, 0.3] & 1 & 19& [0.86, 0.07, 0.07] & [0.86, 0.14] & [0.8, 0.2] &[0.8, 0.2] & 3 & 3 &59&6&16&19&0&0&30&1951\\
42&$p(38.30.3.15.16.3.15)$ & 2 & [0.7, 0.3] & 1 & 19& [0.86, 0.07, 0.07] & [0.86, 0.14] & [0.7, 0.3] &[0.7, 0.3] & 3 & 3 &39&5&38&18&0&0&29&1215\\
43&$p(38.30.3.15.16.3.15)$ & 2 & [0.7, 0.3] & 1 & 19& [0.86, 0.07, 0.07] & [0.8, 0.2] & [0.8, 0.2] &[0.8, 0.2] & 3 & 3 &26&5&31&18&16&4&30&35358\\
44&$p(38.30.3.15.16.3.15)$ & 2 & [0.7, 0.3] & 1 & 19& [0.86, 0.07, 0.07] & [0.8, 0.2] & [0.7, 0.3] &[0.7, 0.3] & 3 & 3 &39&5&38&18&0&0&30&20317\\
45&$p(38.30.3.15.16.3.15)$ & 2 & [0.7, 0.3] & 1 & 19& [0.82, 0.09, 0.09]  & [0.86, 0.14] & [0.8, 0.2] &[0.8, 0.2] & 3 & 3 &44&6&28&22&0&0&30&6579\\
46&$p(38.30.3.15.16.3.15)$ & 2 & [0.7, 0.3] & 1 & 19& [0.82, 0.09, 0.09]  & [0.86, 0.14] & [0.7, 0.3] &[0.7, 0.3] & 3 & 3 &42&4&25&12&13&4&30&34954\\
47&$p(38.30.3.15.16.3.15)$ & 2 & [0.7, 0.3] & 1 & 19& [0.82, 0.09, 0.09]  & [0.8, 0.2] & [0.8, 0.2] &[0.8, 0.2] & 3 & 3 &14&4&37&11&20&14&30&38044\\
48&$p(38.30.3.15.16.3.15)$ & 2 & [0.7, 0.3] & 1 & 19& [0.82, 0.09, 0.09]  & [0.8, 0.2] & [0.7, 0.3] &[0.7, 0.3] & 3 & 3 &46&5&27&17&2&3&30&21212\\

\hdashline
49&$p(38.30.3.15.16.4.15)$ & 2 & [0.7, 0.3] & 2 & 19& [0.82, 0.09, 0.09]  & [0.86, 0.14] & [0.8, 0.2] &[0.8, 0.2] & 3 & 3 &40&5&40&15&0&0&23&1589\\
50&$p(38.30.3.15.16.4.15)$ & 2 & [0.7, 0.3] & 2 & 19& [0.82, 0.09, 0.09]  & [0.86, 0.14] & [0.7, 0.3] &[0.7, 0.3] & 3 & 3 &40&3&47&10&0&0&28&1351\\
51&$p(38.30.3.15.16.4.15)$ & 2 & [0.7, 0.3] & 2 & 19& [0.82, 0.09, 0.09]  & [0.8, 0.2] & [0.8, 0.2] &[0.8, 0.2] & 3 & 3 &38&6&37&19&0&0&30&2984\\
52&$p(38.30.3.15.16.4.15)$ & 2 & [0.7, 0.3] & 2 & 19& [0.82, 0.09, 0.09]  & [0.8, 0.2] & [0.7, 0.3] &[0.7, 0.3] & 3 & 3 &32&3&58&7&0&0&23&1248\\
53&$p(38.30.3.15.16.4.15)$ & 2 & [0.7, 0.3] & 2 & 19& [0.78, 0.11, 0.11] & [0.86, 0.14] & [0.8, 0.2] &[0.8, 0.2] & 3 & 3 &53&6&25&16&0&0&25&1862\\
54&$p(38.30.3.15.16.4.15)$ & 2 & [0.7, 0.3] & 2 & 19& [0.78, 0.11, 0.11] & [0.86, 0.14] & [0.7, 0.3] &[0.7, 0.3] & 3 & 3 &51&5&32&12&0&0&24&1883\\
55&$p(38.30.3.15.16.4.15)$ & 2 & [0.7, 0.3] & 2 & 19& [0.78, 0.11, 0.11] & [0.8, 0.2] & [0.8, 0.2] &[0.8, 0.2] & 3 & 3 &40&5&41&14&0&0&26&2218\\
56&$p(38.30.3.15.16.4.15)$ & 2 & [0.7, 0.3] & 2 & 19& [0.78, 0.11, 0.11] & [0.8, 0.2] & [0.7, 0.3] &[0.7, 0.3] & 3 & 3 &32&5&49&14&0&0&24&2382\\
57&$p(38.30.3.15.16.4.15)$ & 2 & [0.7, 0.3] & 1 & 19& [0.86, 0.07, 0.07] & [0.86, 0.14] & [0.8, 0.2] &[0.8, 0.2] & 3 & 3 &54&5&26&15&0&0&26&2028\\
58&$p(38.30.3.15.16.4.15)$ & 2 & [0.7, 0.3] & 1 & 19& [0.86, 0.07, 0.07] & [0.86, 0.14] & [0.7, 0.3] &[0.7, 0.3] & 3 & 3 &57&5&26&12&0&0&30&2493\\
59&$p(38.30.3.15.16.4.15)$ & 2 & [0.7, 0.3] & 1 & 19& [0.86, 0.07, 0.07] & [0.8, 0.2] & [0.8, 0.2] &[0.8, 0.2] & 3 & 3 &53&5&29&13&0&0&30&18086\\
60&$p(38.30.3.15.16.4.15)$ & 2 & [0.7, 0.3] & 1 & 19& [0.86, 0.07, 0.07] & [0.8, 0.2] & [0.7, 0.3] &[0.7, 0.3] & 3 & 3 &15&4&31&14&22&14&30&37324\\
61&$p(38.30.3.15.16.4.15)$ & 2 & [0.7, 0.3] & 1 & 19& [0.82, 0.09, 0.09]  & [0.86, 0.14] & [0.8, 0.2] &[0.8, 0.2] & 3 & 3 &22&5&36&18&11&8&30&31803\\
62&$p(38.30.3.15.16.4.15)$ & 2 & [0.7, 0.3] & 1 & 19& [0.82, 0.09, 0.09]  & [0.86, 0.14] & [0.7, 0.3] &[0.7, 0.3] & 3 & 3 &15&4&30&10&24&17&30&39081\\
63&$p(38.30.3.15.16.4.15)$ & 2 & [0.7, 0.3] & 1 & 19& [0.82, 0.09, 0.09]  & [0.8, 0.2] & [0.8, 0.2] &[0.8, 0.2] & 3 & 3 &8&4&23&11&31&23&30&39821\\
64&$p(38.30.3.15.16.4.15)$ & 2 & [0.7, 0.3] & 1 & 19& [0.82, 0.09, 0.09]  & [0.8, 0.2] & [0.7, 0.3] &[0.7, 0.3] & 3 & 3 &11&4&19&14&36&16&30&38628\\

    \hline
    \end{tabular}
}
    \label{tab:my_label}
\end{table}

\begin{table}[]
    \caption{Computational experiments - large instances ($n_i = 50, \; n_j = 40$)}\label{large}
    \scalebox{0.8}{
    \centering
    \begin{tabular}{|c|c|cccccccccc|cccccccc|}
    \hline
    \multicolumn{2}{|c|}{{Instance}}&\multicolumn{10}{c|}{{Data}}&\multicolumn{8}{c|}{{Output}}\\
    \hline
    N&$p(n_i.n_j.n_t.n_k.n_\ell. n_p.n_q)$&$d$&$u_i$&$e_{ijt}$&$c_i$&$l_{ik\ell}$&$m_{k\ell p}$&$v_{i\ell}$&$h_{i\ell}$&$r_{iq}$&$t_{iq}$&$|N|$&$|I|$&$skip^N$&$skip^I$&$time^N$&$time^I$&$g$&CPU\\
    \hline
65&$p(50.40.3.15.16.3.15)$ & 2 & [0.7, 0.3] & 2 & 25& [0.82, 0.09, 0.09]  & [0.86, 0.14] & [0.8, 0.2] &[0.8, 0.2] & 3 & 3 &28&5&52&15&0&0&28&1366\\
66&$p(50.40.3.15.16.3.15)$ & 2 & [0.7, 0.3] & 2 & 25& [0.82, 0.09, 0.09]  & [0.86, 0.14] & [0.7, 0.3] &[0.7, 0.3] & 3 & 3 &51&4&36&9&0&0&32&2633\\
67&$p(50.40.3.15.16.3.15)$ & 2 & [0.7, 0.3] & 2 & 25& [0.82, 0.09, 0.09]  & [0.8, 0.2] & [0.8, 0.2] &[0.8, 0.2] & 3 & 3 &40&4&43&13&0&0&34&2744\\
68&$p(50.40.3.15.16.3.15)$ & 2 & [0.7, 0.3] & 2 & 25& [0.82, 0.09, 0.09]  & [0.8, 0.2] & [0.7, 0.3] &[0.7, 0.3] & 3 & 3 &36&5&42&17&0&0&35&6862\\
69&$p(50.40.3.15.16.3.15)$ & 2 & [0.7, 0.3] & 2 & 25& [0.78, 0.11, 0.11] & [0.86, 0.14] & [0.8, 0.2] &[0.8, 0.2] & 3 & 3 &39&5&43&11&1&1&35&22050\\
70&$p(50.40.3.15.16.3.15)$ & 2 & [0.7, 0.3] & 2 & 25& [0.78, 0.11, 0.11] & [0.86, 0.14] & [0.7, 0.3] &[0.7, 0.3] & 3 & 3 &49&5&30&16&0&0&32&2719\\
71&$p(50.40.3.15.16.3.15)$ & 2 & [0.7, 0.3] & 2 & 25& [0.78, 0.11, 0.11] & [0.8, 0.2] & [0.8, 0.2] &[0.8, 0.2] & 3 & 3 &40&3&45&9&0&3&33&29643\\
72&$p(50.40.3.15.16.3.15)$ & 2 & [0.7, 0.3] & 2 & 25& [0.78, 0.11, 0.11] & [0.8, 0.2] & [0.7, 0.3] &[0.7, 0.3] & 3 & 3 &40&5&39&15&1&0&37&21446\\
73&$p(50.40.3.15.16.3.15)$ & 2 & [0.7, 0.3] & 1 & 25& [0.86, 0.07, 0.07] & [0.86, 0.14] & [0.8, 0.2] &[0.8, 0.2] & 3 & 3 &48&5&24&19&4&0&40&18000\\
74&$p(50.40.3.15.16.3.15)$ & 2 & [0.7, 0.3] & 1 & 25& [0.86, 0.07, 0.07] & [0.86, 0.14] & [0.7, 0.3] &[0.7, 0.3] & 3 & 3 &39&5&39&17&0&0&36&13145\\
75&$p(50.40.3.15.16.3.15)$ & 2 & [0.7, 0.3] & 1 & 25& [0.86, 0.07, 0.07] & [0.8, 0.2] & [0.8, 0.2] &[0.8, 0.2] & 3 & 3 &17&5&42&14&14&8&40&34408\\
76&$p(50.40.3.15.16.3.15)$ & 2 & [0.7, 0.3] & 1 & 25& [0.86, 0.07, 0.07] & [0.8, 0.2] & [0.7, 0.3] &[0.7, 0.3] & 3 & 3 &42&5&35&16&0&2&40&20341\\
77&$p(50.40.3.15.16.3.15)$ & 2 & [0.7, 0.3] & 1 & 25& [0.82, 0.09, 0.09]  & [0.86, 0.14] & [0.8, 0.2] &[0.8, 0.2] & 3 & 3 &21&4&40&13&16&6&40&35309\\
78&$p(50.40.3.15.16.3.15)$ & 2 & [0.7, 0.3] & 1 & 25& [0.82, 0.09, 0.09]  & [0.86, 0.14] & [0.7, 0.3] &[0.7, 0.3] & 3 & 3 &34&5&32&14&9&6&40&31078\\
79&$p(50.40.3.15.16.3.15)$ & 2 & [0.7, 0.3] & 1 & 25& [0.82, 0.09, 0.09]  & [0.8, 0.2] & [0.8, 0.2] &[0.8, 0.2] & 3 & 3 &15&5&20&15&21&24&39&37261\\
80&$p(50.40.3.15.16.3.15)$ & 2 & [0.7, 0.3] & 1 & 25& [0.82, 0.09, 0.09]  & [0.8, 0.2] & [0.7, 0.3] &[0.7, 0.3] & 3 & 3 &3&3&0&11&39&44&40&41412\\

\hdashline
81&$p(50.40.3.15.16.4.15)$ & 2 & [0.7, 0.3] & 2 & 25& [0.82, 0.09, 0.09]  & [0.86, 0.14] & [0.8, 0.2] &[0.8, 0.2] & 3 & 3 &44&5&37&14&0&0&31&2840\\
82&$p(50.40.3.15.16.4.15)$ & 2 & [0.7, 0.3] & 2 & 25& [0.82, 0.09, 0.09]  & [0.86, 0.14] & [0.7, 0.3] &[0.7, 0.3] & 3 & 3 &39&6&34&21&0&0&26&2726\\
83&$p(50.40.3.15.16.4.15)$ & 2 & [0.7, 0.3] & 2 & 25& [0.82, 0.09, 0.09]  & [0.8, 0.2] & [0.8, 0.2] &[0.8, 0.2] & 3 & 3 &38&4&48&10&0&0&34&9626\\
84&$p(50.40.3.15.16.4.15)$ & 2 & [0.7, 0.3] & 2 & 25& [0.82, 0.09, 0.09]  & [0.8, 0.2] & [0.7, 0.3] &[0.7, 0.3] & 3 & 3 &44&5&35&16&0&0&32&3543\\
85&$p(50.40.3.15.16.4.15)$ & 2 & [0.7, 0.3] & 2 & 25& [0.78, 0.11, 0.11] & [0.86, 0.14] & [0.8, 0.2] &[0.8, 0.2] & 3 & 3 &43&4&42&11&0&0&35&8356\\
86&$p(50.40.3.15.16.4.15)$ & 2 & [0.7, 0.3] & 2 & 25& [0.78, 0.11, 0.11] & [0.86, 0.14] & [0.7, 0.3] &[0.7, 0.3] & 3 & 3 &32&5&50&13&0&0&32&2459\\
87&$p(50.40.3.15.16.4.15)$ & 2 & [0.7, 0.3] & 2 & 25& [0.78, 0.11, 0.11] & [0.8, 0.2] & [0.8, 0.2] &[0.8, 0.2] & 3 & 3 &49&5&31&15&0&0&34&20443\\
88&$p(50.40.3.15.16.4.15)$ & 2 & [0.7, 0.3] & 2 & 25& [0.78, 0.11, 0.11] & [0.8, 0.2] & [0.7, 0.3] &[0.7, 0.3] & 3 & 3 &40&4&41&11&2&2&39&31311\\
89&$p(50.40.3.15.16.4.15)$ & 2 & [0.7, 0.3] & 1 & 25& [0.86, 0.07, 0.07] & [0.86, 0.14] & [0.8, 0.2] &[0.8, 0.2] & 3 & 3 &12&4&12&14&33&25&40&39083\\
90&$p(50.40.3.15.16.4.15)$ & 2 & [0.7, 0.3] & 1 & 25& [0.86, 0.07, 0.07] & [0.86, 0.14] & [0.7, 0.3] &[0.7, 0.3] & 3 & 3 &34&5&29&16&10&6&36&30877\\
91&$p(50.40.3.15.16.4.15)$ & 2 & [0.7, 0.3] & 1 & 25& [0.86, 0.07, 0.07] & [0.8, 0.2] & [0.8, 0.2] &[0.8, 0.2] & 3 & 3 &10&4&18&10&17&41&40&40127\\
92&$p(50.40.3.15.16.4.15)$ & 2 & [0.7, 0.3] & 1 & 25& [0.86, 0.07, 0.07] & [0.8, 0.2] & [0.7, 0.3] &[0.7, 0.3] & 3 & 3 &16&4&37&11&10&22&39&38514\\
93&$p(50.40.3.15.16.4.15)$ & 2 & [0.7, 0.3] & 1 & 25& [0.82, 0.09, 0.09]  & [0.86, 0.14] & [0.8, 0.2] &[0.8, 0.2] & 3 & 3 &11&5&22&14&16&32&40&38139\\
94&$p(50.40.3.15.16.4.15)$ & 2 & [0.7, 0.3] & 1 & 25& [0.82, 0.09, 0.09]  & [0.86, 0.14] & [0.7, 0.3] &[0.7, 0.3] & 3 & 3 &19&5&40&13&13&10&40&35060\\
95&$p(50.40.3.15.16.4.15)$ & 2 & [0.7, 0.3] & 1 & 25& [0.82, 0.09, 0.09]  & [0.8, 0.2] & [0.8, 0.2] &[0.8, 0.2] & 3 & 3 &7&3&17&9&31&33&40&41043\\
96&$p(50.40.3.15.16.4.15)$ & 2 & [0.7, 0.3] & 1 & 25& [0.82, 0.09, 0.09]  & [0.8, 0.2] & [0.7, 0.3] &[0.7, 0.3] & 3 & 3 &6&3&16&10&23&42&40&40999\\

    \hline
    \end{tabular}
}

\end{table}

\end{landscape}

\end{document}